\documentclass{article}

\usepackage{graphicx}
\usepackage{verbatim}
\usepackage{enumitem}
\usepackage{caption}

\usepackage{float}
\usepackage{fancybox}
\usepackage{textcomp}
\usepackage{color}

\usepackage{amsmath}
\usepackage{amsfonts}
\usepackage{amssymb}
\usepackage{amsthm}
\usepackage{makeidx,bbm}

\usepackage{mathrsfs}
\usepackage{eurosym}
\usepackage[T1]{fontenc}

\usepackage{makeidx}

\usepackage[utf8]{inputenc}
\usepackage[english]{babel}
\usepackage{young}
\usepackage[vcentermath]{youngtab}
\usepackage{tikz}
\usetikzlibrary{calc}
\usepackage[nomessages]{fp}
\usepackage{amsthm}
\usepackage{adjustbox}
\usepackage{amsmath}
\usepackage[affil-it]{authblk}
\usepackage{mathtools}
\usepackage{chemfig}
\usepackage{parallel,enumitem}
\usepackage{pdfcolparallel}
\usepackage{setspace}
\usepackage{fancyhdr}

\usepackage{romannum}

\theoremstyle{plain}
\newtheorem{thm}{Theorem}[section]

\newtheorem{cor}{Corollary}[thm]
\newtheorem{conj}[thm]{Conjecture}

\newtheorem{remark}[thm]{Remark}
\theoremstyle{definition}
\newtheorem{defn}[thm]{Definition}
\newtheorem{exmp}[thm]{Example}

\newcommand{\R}{\mathbb{R}}
\newcommand{\vs}{\vec w}

\evensidemargin \oddsidemargin
\newcommand*{\Scale}[2][4]{\scalebox{#1}{$#2$}}%

\begin{document} 

\title{A new approach to $e$-positivity for \\ Stanley's chromatic functions} 
\author{Alexander Paunov, Andr\'as Szenes} 
\date{February 2016} 
\maketitle

\pagestyle{plain}

\begin{abstract}

In this paper, we study positivity phenomena
  for the $e$-coefficients of Stanley's chromatic function of a graph.
  We introduce a new combinatorial object: the {\em correct} sequences of unit
interval orders, and using these, in certain cases,  we succeed to
construct combinatorial models of
the coefficients appearing in Stanley's conjecture. Our main result is
the proof of positivity of the
  coefficients $c_{n-k,1^k}$, $c_{n-2,2}$, $c_{n-3,2,1}$ and
  $c_{2^k,1^{n-2k}}$ of the expansion of the chromatic symmetric
  function in terms of the basis of the elementary symmetric
  polynomials for the case of $(3+1)$-free posets.  

\end{abstract}

\section{Introduction}

Let $G$ be a finite graph, $V(G)$ - the set of vertices of $G$, $E(G)$ - the set of edges of $G$. 

\begin{defn} \label{coloring} A {\em proper coloring} $c$ of  $G$ is
  a map $$c:V\rightarrow\mathbb{N}$$ such that no two adjacent
  vertices are colored in the same color.
\end{defn}

For each coloring $c$ we define a monomial $$x^c = \prod_{v\in
  V}x_{c(v)},$$ where $x_1, x_2, ..., x_n,...$ are commuting
variables.  We denote by $\Pi(G)$ the set of all proper colorings
of $G$, and by  $\Lambda$ the ring of symmetric functions in the infinite set
of variables $\{x_1, x_2,...\}.$

In \cite{Stanley95a}, Stanley defined the chromatic symmetric function of a graph.

\begin{defn} \label{chromfunction} The \em chromatic symmetric
  function \normalfont $X_G\in\Lambda$ of a graph $G$ is the sum of the monomials $x^c$
  over all proper colorings of
  $G$: $$X_G=\sum\limits_{c\in\Pi(G)}x^c.$$
\end{defn}

\begin{defn} \label{efunc} Denote by $e_m$ the $m$-th elementary
  symmetric function:
  $$e_m = \sum\limits_{i_1<i_2<...<i_m}x_{i_1}\cdot
  x_{i_2}\cdot...\cdot x_{i_m},$$
  where $i_1,..,i_k\in \mathbb{N}$.  Given a non-increasing sequence of
  positive integers (we will call these {\em partitions})
  $$\lambda = (\lambda_1\geq \lambda_2\geq...\geq\lambda_k),\ \lambda_i\in
  \mathbb{N},$$
  we define the elementary symmetric function
  $e_{\lambda} = \prod\limits_{i=1}^k e_{\lambda_i}.$ These functions
  form a basis of $\Lambda.$
\end{defn}
For a natural number $k$, we denote by $1^k$ the partition $\lambda$ of length $k$, where $$\lambda_1=\lambda_2=...=\lambda_k=1.$$

\begin{defn} \label{epos} A symmetric function $X\in \Lambda$ is \em
  $e$-positive \normalfont if it has non-negative coefficients in the
  basis of the elementary symmetric functions.
\end{defn}

\begin{defn} \label{pfunc} Denote by $p_m$ the $m$-th power sum
  symmetric function: $$p_m =
  \sum\limits_{i\in\mathbb{N}}x^m_{i}.$$ Given a
  partition $\lambda = (\lambda_1\geq \lambda_2\geq...\geq\lambda_k)$, we define the power sum
  symmetric function
  $p_{\lambda} = \prod\limits_{i=1}^k p_{\lambda_i}.$ These functions also
  form a basis of $\Lambda.$
\end{defn}

\begin{defn} \label{mfunc}
Given a partition $\lambda = (\lambda_1\geq \lambda_2\geq...\geq\lambda_k)$,  we define the monomial symmetric function $$m_\lambda=\sum\limits_{i_1<i_2<...<i_k}\sum\limits_{\lambda'\in S_k(\lambda)}x_{i_1}^{\lambda_{1}'}\cdot
  x_{i_2}^{\lambda_{1}'}\cdot...\cdot x_{i_k}^{\lambda_{k}'},$$
where the inner sum is taken over the set of all permutations of the sequence $\lambda$, denoted by $S_k(\lambda)$.
\end{defn}

\begin{exmp}
  The chromatic symmetric function of $K_n$, the complete graph on $n$
  vertices, is $e$-positive: $X_{K_n} = n!\,e_n$.
\end{exmp}

\begin{defn} \label{incgraph} For a poset $P$, the \em incomparability
  graph\normalfont, $\textnormal{inc}(P)$, is the graph with elements
  of $P$ as vertices, where two vertices are connected if and only if
  they are not comparable in $P$.
\end{defn} 

\begin{defn} \label{nplusmfree} Given a pair of natural numbers
  $a,b\in\mathbb{N}^2$, we say that a poset $P$ is \em (a+b)-free
  \normalfont if it does not contain a length-$a$ and a length-$b$
  chain, whose elements are  incomparable.
\end{defn} 

\begin{defn} A unit interval order (UIO) is a partially ordered set
  which is isomorphic to a finite subset of $U\subset\R$ with the following poset structure:
\[ \text{for } u,w\in U:\    u\succ w \text{ iff } u\ge w+1.
\]
Thus $u$ and $w$ are incomparable precisely when $|u-w|<1$ and we will
use the notation $u\sim w$ in this case. 
\end{defn}
\begin{thm}[Scott-Suppes \cite{Scott-Suppes54}]\label{S_S}
A finite poset $P$ is a UIO if and only if it is $(2+2)$- and $(3+1)$-free.
\end{thm}

Stanley~\cite{Stanley95a} initiated the study of incomparability
graphs of $(3+1)$-free partially ordered sets. Analyzing the chromatic
symmetric functions of these incomparability graphs,
Stanley~\cite{Stanley95a} stated the following positivity conjecture.

\begin{conj}[Stanley] \label{eposconj}
If $P$ is a $(3+1)$-free poset, then $X_{\textnormal{inc}(P)}$ is $e$-positive.
\end{conj}

For a graph $G$ let us denote by ${c_\lambda}(G)$ the coefficients of $X_G$ with respect to the $e$-basis. We omit the index $G$ whenever this causes no confusion:
$$X_G=\sum\limits_{\lambda}c_{\lambda}e_\lambda.$$

Conjecture~\ref{eposconj} has been verified with the help of computers
for up to 20-element posets~\cite{Guay-Paquet13}. In 2013,
Guay-Paquet~\cite{Guay-Paquet13} showed that to prove this conjecture,
it would be sufficient to verify it for the case of $(3+1)$- and
$(2+2)$-free posets, i.e. for unit interval orders (see Theorem~\ref{S_S}). More precisely:

\begin{thm}[Guay-Paquet]\label{G_P}
  Let $P$ be a $(3+1)$-free poset. Then, $X_\mathrm{inc}(P)$ is a
  convex combination of the chromatic symmetric
  functions $$\{X_\mathrm{inc}(P')\ |\ P'\ \mathrm{is}\ \mathrm{a}\
  (3+1)\mathrm{-}\ \mathrm{and}\ (2+2)\mathrm{-free}\ \mathrm{poset}
  \}.$$

\end{thm}

The strongest general result in this direction is that of Gasharov
\cite{Gasharov94}.

\begin{defn} \label{sfunc} For a partition $\lambda = (\lambda_1\geq
  \lambda_2\geq...\geq\lambda_k)$, define the Schur
  functions \em $s_{\lambda}=\mathrm{det}(e_{\lambda_i^*+j-i})_{i,j}$,
  \normalfont where $\lambda^*$ is the conjugate partition to
  $\lambda$. The functions $\{ s_{\lambda}\}$ form a basis of
  $\Lambda$.
\end{defn}

\begin{defn} \label{spos}
A symmetric polynomial $X$ is \em $s$-positive \normalfont if it has non-negative coefficients in the basis of Schur functions.
\end{defn}

Obviously, a product of $e$-positive functions is $e$-positive. This also holds for $s$-positive functions. Thus, the equality $e_n=s_{1^n}$ implies that $e$-positive functions are $s$-positive, and thus $s$-positivity is weaker than $e$-positivity.

\begin{thm}[Gasharov] \label{sposthm} 
If $P$ is a $(3+1)$-free poset, then $X_{\textnormal{inc}(P)}$ is $s$-positive.
\end{thm}

Gasharov proved $s$-positivity by  constructing so-called $P$-tableau and finding a one-to-one correspondence between these tableau and $s$-coefficients~\cite{Gasharov94}. However, $e$-positivity conjecture~\ref{epos} is still open. 
The strongest known result on the $e$-coefficients was obtained by Stanley in~\cite{Stanley95a}. He showed that sums of $e$-coefficients over the partitions of fixed length are non-negative:

\begin{thm}[Stanley]
For a finite graph $G$ and $j\in\mathbb{N}$, suppose $$X_G=\sum\limits_{\lambda}c_{\lambda}e_\lambda,$$
and let $\text{sink}(G,j)$ be the number of acyclic orientation of $G$ with $j$ sinks. Then
$$\text{sink}(G,j)=\sum\limits_{l(\lambda)=j}c_{\lambda}.$$
\end{thm}

\begin{remark}
By taking $j=1$, it follows from the theorem that $c_n$ is non-negative.
\end{remark}

Stanley in~\cite{Stanley95a} showed that for $n\in\mathbb{N}$ and the unit interval order $P_n=\{\frac{i}{2}\}_{i=1}^n$, the corresponding $X_{\text{inc}(P_n)}$ is $e$-positive, while $e$-positivity for the UIOs $$P_{n,k}=\bigg\{\frac{i}{k+1}\bigg\}_{i=1}^n$$ with $k>1$ has not yet been proven. It was checked for small $n$ and some $k$ (see~\cite{Stanley95a}).


Next, we introduce {\em correct sequences} (abbreviated as {\em corrects}),
defined below. These play a major role in the article.

\begin{defn}
Let U be a UIO.  We will call a sequence $\vs = (w_1,\dots, w_k)$ of elements of $U$ {\em
   correct} if
 \begin{itemize}
 \item 
$w_i\not\succ w_{i+1}$ for $i=1,2,\dots,k-1$ 
\item and for
 each $j=2,\dots,k$, there exists $i<j$ such that $w_i\not\prec w_j$.
 \end{itemize}
\end{defn}
Every sequence of length 1 is correct, and sequence $(w_1,w_2)$ is
correct precisely when $w_1\sim w_2$. The second condition (supposing that the first one holds) may be reformulated as follows:
for each $j=1,\dots k$, the subset $\{w_1,\dots,w_j\}\subset U$ is connected
with respect to the graph structure~${(U,\sim)}$.
Using this notation, we prove the following theorems.
\begin{thm}\label{eposn} 
Let $X_{\text{inc}(U)}=\sum\limits_{\lambda}c_\lambda e_\lambda$ be a chromatic symmetric function of the $n$-element unit interval order $U$. Then $c_n$ is equal to the number of corrects of length $n$, in which every element of $U$ is used exactly once.
\end{thm}
\begin{cor}
Let $X_{\text{inc}(P)}=\sum\limits_{\lambda}c_\lambda e_\lambda$ be a chromatic symmetric function of  $n$-element $(3+1)$-free poset $P$, then $c_n$ is a nonnegative integer.
\end{cor}
Indeed, positivity for the general case follows from Theorem~\ref{G_P}, which presents the chromatic symmetric function of a $(3+1)$-free poset as a convex combination of the chromatic symmetric functions of unit interval orders.

Stanley~\cite{Stanley95b} and Chow~\cite{Chow95} showed the positivity of $c_n$ for $(3+1)$-free posets using combinatorial techniques, and linked $e$-coefficients with the acyclic orientations of the incomparability graphs. The construction of corrects not only serves this purpose for UIOs (see~\cite{Paunov16b}), but also creates a new approach, which allows us to obtain the following new result:

\begin{thm} \label{eposn21} 
Let $X_{\text{inc}(P)}=\sum\limits_{\lambda}c_\lambda e_\lambda$ be a chromatic symmetric function of the $(3+1)$-free poset~$P$, and $k\in\mathbb{N}$. Then $c_{n-k,1^k}$, $c_{n-2,2}$, $c_{n-3,2,1}$ and $c_{2^k,1^{n-2k}}$ are non-negative integers.
\end{thm}

The proofs of Theorem~\ref{eposn} and  Theorem~\ref{sposthm},  and positivity of correspondent $G$-power sum symmetric functions and  Schur $G$-symmetric functions can be found in~\cite{Paunov16} and~~\cite{Paunov16b}. The article is structured as follows: in Section~\ref{Ghom}, we describe
the $G$-homomorphism introduced by Stanley in~\cite{Stanley95b}, which
is essential for our approach. Positivity of $c_{n-k,1^k}$, $c_{n-2,2}$, $c_{n-3,2,1}$ and $c_{2^k,1^{n-2k}}$ (Theorem~\ref{eposn21}) is proven in Section~3.

\textbf {Acknowledgements.}  We are grateful to
Emanuele Delucchi and Bart Vandereycken for their help and useful
discussions.

\section{Stanley's $G$-homomorphism}\label{Ghom}
For a graph $G$, Stanley \cite[p.~6]{Stanley95b} defined $G$-analogues of the standard families of symmetric functions. Let $G$ be a finite graph with vertex set $V(G)=\{v_1,...,v_n\}$ and edge set $E(G)$. We will think of the elements of $V(G)$ as commuting variables.

\begin{defn} \label{eG}
For a positive integer $i$, $1\leq i, \leq n$, we define the \em $G$-analogues \normalfont of the elementary symmetric polynomials, or \em the elementary $G$-symmetric polynomials\normalfont, as follows
 $$e_i^G =\sum\limits_{\substack{\#S=i\\
     S-\mathrm{stable}}}\prod\limits_{v\in S}v,$$ where the sum is
 taken over all $i$-element subsets $S$ of $V$, in which no two
 vertices form an edge, i.e. stable subsets. We set $e_0^G=1$, and $e_i^G=0$ for $i<0$.
\end{defn}
Note that these polynomials are not necessarily symmetric.

Let $\Lambda_G\subset\mathbb{R}[v_1,...,v_n]$ be the subring generated
by $\{e_i^G\}_{i=1}^{n}$.  The map $e_i\mapsto e_i^G$ extends to a
ring homomorphism $\phi_G: \Lambda\rightarrow\Lambda_G$, called the
{\em $G$-homomorphism}. For $f\in \Lambda$, we will use the notation $f^G$ for $\phi_G(f)$.

\begin{exmp}
Given a partition $\lambda = \lambda_1\geq \lambda_2\geq...\geq\lambda_k,\ k\in \mathbb{N},$ we have $$e_{\lambda}^G = \prod\limits_{i=1}^k e_i^G,$$
$$s_{\lambda}^G=\mathrm{det}(e_{\lambda_i^*+j-i}^G).$$
\end{exmp}

For an integer function $\alpha: V\rightarrow \mathbb{N}$ and $f^G\in\Lambda_G$, let
$$v^\alpha = \prod\limits_{v\in V}v^{\alpha(v)},$$ 
 and $[v^\alpha]f^G$ stands for the coefficient of $v^\alpha$ in the polynomial $f^G\in\Lambda_G$.

Let $G^\alpha$ denote the graph, obtained by replacing every vertex $v$ of $G$ by the complete subgraph of size $\alpha(v)$: $K_{\alpha(v)}^v$. Given vertices $u$ and $v$ of $G$, a vertex of $K_{\alpha(v)}^v$ is connected to a vertex of $K_{\alpha(u)}^u$ if and only if $u$ and $v$ form an edge in $G$.

Considering the Cauchy product \cite[ch.~4.2]{Macdonald79}, Stanley \cite[p.~6]{Stanley95b} found a connection between the $G$-analogues of symmetric functions and $X_G$. Following Stanley~\cite{Stanley95b}, we set
$$T(x,v) = \sum\limits_\lambda m_\lambda(x)e^G_\lambda(v),$$
where the sum is taken over all partitions. Then 

\begin{equation}\label{gnechrom}
[v^\alpha]T(x,v)\prod\limits_{v\in V}\alpha(v)! =X_{G^\alpha}.
\end{equation}
Using the Cauchy identity 
$$\sum\limits_\lambda s_\lambda(x)s_{\lambda^*}(y)=\sum\limits_\lambda m_\lambda(x)e_\lambda(y) = \sum\limits_\lambda e_\lambda(x)m_\lambda(y)$$
and applying the $G$-homomorphism, one obtains:
\begin{equation}\label{GCauchy}
T(x,v) = \sum\limits_\lambda m_\lambda(x)e^G_\lambda(v) = \sum\limits_\lambda s_\lambda(x)s^G_{\lambda^*}(v)=T(v,x) = \sum\limits_\lambda e_\lambda(x)m^G_\lambda(v).
\end{equation}

An immediate consequence of the formulas \eqref{gnechrom} and
\eqref{GCauchy} is the following result of Stanley:
\begin{thm}[Stanley]\label{poscrit}
 For every finite graph G
\begin{enumerate}
\item $X_{G^\alpha}$ is s-positive for every
  $\alpha:V(G)\rightarrow\mathbb{N}$ if and only if $s_\lambda^G\in
  \mathbb{N}[V(G)]$ for every partition $\lambda$.
\item $X_{G^\alpha}$ is e-positive for every $\alpha:V(G)\rightarrow\mathbb{N}$ if and only if $m_\lambda^G\in \mathbb{N}[V(G)]$ for every partition $\lambda$. 
\end{enumerate}
\end{thm}

\begin{remark}\label{c_m}
If $X_{G^\alpha}=\sum\limits_{\lambda}c^\alpha_{\lambda}e_\lambda,$ then $c_\lambda^\alpha=[v^\alpha]m^G_\lambda.$ Hence, monomial positivity of $m^G_\lambda$ is equivalent to the positivity of $c_\lambda^\alpha$ for every $\alpha$.
\end{remark}

The proofs of positivity of $G$-power sum symmetric functions and  Schur $G$-symmetric functions for the case of unit interval orders can be found in~\cite{Paunov16}.

\section{Proofs of the theorems}\label{proofs}
It follows from Theorem \ref{poscrit} that to prove that the graph $G$ is $e$-positive, it is enough to show the monomial positivity of its monomial $G$-symmetric functions. On the other hand, Guay-Paquet in Theorem~\ref{G_P} showed that it is sufficient to check $e$-positivity for unit interval orders, in order to prove it for the general case of $(3+1)$-free posets. Therefore, in the following section~\ref{proofs} we analyze the functions $m_{\lambda}^G$ for the case $G=\text{inc}(U),$ where $U$ is UIO.

Let us repeat the definition of a central notion for our work, that of
correct sequences of elements of a unit interval order.
\begin{defn}
 Let $(U,\prec)$ be a unit interval order, and $G=\text{inc}(U)$. We will call a sequence $\vec{w} = (w_1,\dots, w_k)$ of elements of $U$ {\em
   correct} if
 \begin{itemize}
 \item 
$w_i\not\succ w_{i+1}$ for $i=1,2,\dots,k-1$ 
\item and for
 each $j=2,\dots,k$, there exists $i<j$ such that $w_i\not\prec w_j$.
 \end{itemize}
\end{defn}
 
We denote by $P^U_k$ the set of all correct sequences (abbreviated as  {\em
   corrects}) of length $k$. Since $G$ is uniquely defined by $U$, and we are working only with UIO, here and below we use the $U$-index instead of $G$. The $U$-analogues of symmetric functions will be analyzed.

\begin{thm}\label{Ppos}
Let $U$ be a unit interval order and $p_k^U$ the Stanley power-sum function of the corresponding incomparability graph.
Then, for
every natural $k$, we have $$p_k^U=\sum\limits_{\vec{w}\in P^U_k} w_1\cdot...\cdot w_k\ \in N[U],$$ where the sum is taken over all corrects of length $k$.
\end{thm}

The proof of this theorem can be found in~\cite{Paunov16}.


\newcommand{\I}{|}

\newcommand{\w}{w}  
\newcommand{\wi}{y}     
\newcommand{\wii}{u}   

\newcommand{\kap}{q}   
\newcommand{\xix}{\xi}

\newcommand{\ind}{l}

\newcommand{\cph}{\color{red}}   
\newcommand{\cps}{\color{orange}}   

Below, we prove positivity of $m^U_{\ind,1^k}$, $m^U_{\ind,2}$, $m^U_{\ind,2,1},$ and $m^U_{2^\ind,1^k}$. We need the following mild technical generalization of correct
sequences: let $\lambda=(\lambda_1\ge\dots\ge \lambda_k)$ be a
partition of $|\lambda|=\sum\limits_{i=1}^k\lambda_i$. Then, we will call sequence $(w_1,\dots, w_{|\lambda|})$
{\em $\lambda$-correct} if each of the subsequences $(w_1,\dots w_{\lambda_1})$,  $(w_{\lambda_1+1},\dots,
w_{\lambda_1+\lambda_2})$,$\dots$ $(w_{|\lambda|-\lambda_k+1},\dots, w_{|\lambda|})$ are
correct. Introduce the set 
\[   P_\lambda^U=\{\vs=(w_1,\dots w_{\lambda_1}\I w_{\lambda_1+1},\dots,
w_{\lambda_1+\lambda_2}\I\dots\I w_{|\lambda|-\lambda_k+1},\dots, w_{|\lambda|})\ |\ \vs \text{ is
}\lambda\text{-correct }\}
\]
of $\lambda$-correct sequences of length-$|\lambda|$. In particular, $P_\ind^U$ is the set of
$\ind$-corrects of $U$. This definition is consistent with Theorem~\ref{Ppos}, and we have:
$$p_\lambda^U=\prod_{i=1}^kp_{\lambda_i}^U=\sum_{\vec{w}\in P_\lambda^U}w_1\cdot ...\cdot w_{|\lambda|}.$$
For $\vec{w}=(w_1,\dots w_{\ind})\in P^U_\ind$ and $z\in U$ we write $z\succ\vec{w}$, if $z\succ w_i$ for every $1\leq i\leq\ind.$

\begin{thm}\label{Thn1}
Let $$M_{\ind,1}^U = \{(\vs\ |z)\in P_{\ind.1}|\; z\succ\vs\vee z\prec w_\ind\},$$

then $$m^U_{\ind,1}=\sum\limits_{(\vs;z)\in M^U_{\ind,1}}w_1\cdot...\cdot w_\ind\cdot z.$$
\end{thm}
\begin{remark}
According to Remark~\ref{c_m}, this implies $c_{n-1,1}(U)\geq 0$.

\end{remark}

\begin{proof}
Since $P^U_{\ind+1}\subset P^U_{\ind}\times P^U_{1}$, using the following relation
$$m^U_{\ind,1}=p^U_\ind\cdot p^U_1 - p^U_{\ind+1},$$
we have $$P^U_{\ind}\times P^U_{1}\setminus P^U_{\ind+1}=M_{\ind,1}^U,$$
and, as a consequence, $$m^U_{\ind,1}=\sum\limits_{(\vs;z)\in M^U_{\ind,1}}w_1\cdot...\cdot w_\ind\cdot z.$$
\end{proof}

 Next, we introduce the set
\[   E_k^U=\{\vec{\varepsilon}=(\varepsilon_1,\dots, \varepsilon_{k})|\; \varepsilon_i\prec \varepsilon_{i+1}\text{, for }1\leq i<k     \}.
\]

\begin{thm}\label{n1k}
For natural numbers $\ind$ and $k$, let $$M^U_{\ind,1^k}=\{(\vec{w}\ |\vec{\varepsilon})\in P^U_\ind\times E_k^U|\; \varepsilon\prec w_\ind \  \vee\  \varepsilon\succ \vec{w}\text{, for every }  \varepsilon \in \vec{\varepsilon}\},$$
Then,
$$m_{\ind,1^k}^U=\sum_{(\vs;\vec{\varepsilon})\in M^U_{\ind,1^k}}w_1\cdot...\cdot w_{\ind}\cdot \varepsilon_1\cdot...\cdot\varepsilon_k.$$

\end{thm}
\begin{remark}
According to Remark~\ref{c_m}, this implies $c_{n-k,1^k}(U)\geq 0$.

\end{remark}

\begin{proof}

We prove this by induction on $k$.
Note that for $k=1$, the definition of $M^U_{\ind,1^k}$ coincides with $M^U_{\ind,1}$ from Theorem~\ref{Thn1}. Thus, the case $k=1$ follows from Theorem~\ref{Thn1}.

Assume the statement is true for $k$, and consider the standard equation
$$p^U_{\ind}*e^U_{k+1}=m^U_{\ind,1^{k+1}}+m^U_{\ind+1,1^{k}}.$$


Below, we construct a pair of inverse maps, $\phi_{\ind,1^{k+1}}$ and $\psi_{\ind,1^{k+1}}$ from the left part to the right part of the latter equation and vice versa respectively. Every case is followed by a visual illustration.

\bigskip

\begin{minipage}[t]{0.5\textwidth}

{\bf \Romannum{1}. We define 
$$\phi_{\ind,1^{k+1}}:P^U_\ind\times E_{k+1}^U\to M^U_{\ind,1^{k+1}}\sqcup M^U_{\ind+1,1^{k}}$$ as follows:}

Let $(\vec{w}\ |\vec{\varepsilon}\ )\in P^U_\ind\times E_{k+1}^U$
\setstretch{1.75}
\begin{enumerate}
  \item[1.] If $\varepsilon_i\prec w_\ind \  \vee\  \varepsilon_i\succ \vec{w}$ for $1\leq i\leq k+1$, then \setstretch{1.0}
$$\cps \phi_{\ind,1^{k+1}}(\vec{w}\ |\vec{\varepsilon}\ )= (\vec{w}\ ;\vec{\varepsilon}\ )\in M^U_{\ind,1^{k+1}}.$$
\end{enumerate}
\end{minipage}
\vrule \hspace{0.5cm} 
\begin{minipage}[t]{0.5\textwidth}
{\bf  \Romannum{2}. The inverse of map $\psi_{\ind,1^{k+1}}$: \begin{align*}  &\psi^1_{\ind,1^{k+1}}:M^U_{\ind,1^{k+1}}\to P^U_\ind\times E_{k+1}^U;\\ &\psi^2_{\ind,1^{k+1}}:M^U_{\ind+1,1^{k}}\to P^U_\ind\times E_{k+1}^U.  \end{align*} }
Let $ (\vec{w}\ ;\vec{\varepsilon}\ )\in M^U_{\ind,1^{k+1}}$ and $(\vec{w}\ ,z ;\vec{\nu}\ )\in M^U_{\ind+1,1^{k}}$.
\setstretch{1}
\begin{enumerate}
  \item[1.] For $(\vec{w}\ ;\vec{\varepsilon}\ )\in M^U_{\ind,1^{k+1}}$, we have: \setstretch{0.75} $$\cph \psi^1_{\ind,1^{k+1}}(\vec{w}\ ;\vec{\varepsilon}\ )= (\vec{w}\ |\vec{\varepsilon}\ )\in P^U_\ind\times E_{k+1}^U.$$\setstretch{1.0}
\end{enumerate}
\end{minipage}

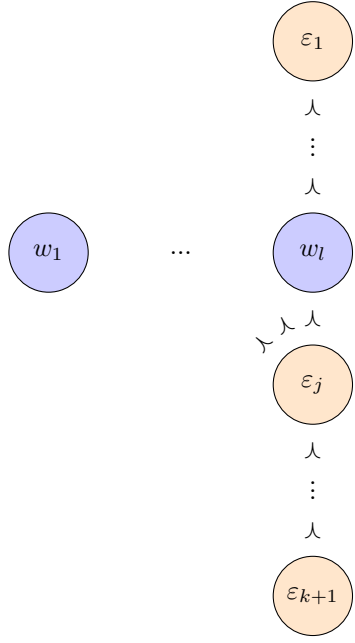
\begin{figure}[H]
\begin{center}
\begin{tikzpicture} [>=latex,every node/.style={minimum width=3em, node distance=4em}]

\node [draw,circle, fill=blue!20] (b) {$\w_1$}; 
\node [node distance=5em,right of=b] (d) {$...$}; 
\node [draw,circle, fill=blue!20, node distance=5em,right of=d] (f) {$\w_{\ind}$}; 
\node [node distance=2.5em, below of=f,rotate=-90] (k) {$\prec$}; 
\node [draw,circle, fill=orange!20, node distance=2.5em, below of=k] (l) {$\varepsilon_j$};
\node [node distance=2.5em, below of=l,rotate=-90] (zl) {$\prec$}; 
\node [node distance=1.5em, below of=zl,rotate=-90] (zl1) {$...$}; 
\node [node distance=1.5em, below of=zl1,rotate=-90] (zl2) {$\prec$}; 
\node [draw,circle, fill=orange!20, node distance=2.5em, below of=zl2] (zq1) {$\varepsilon_{k+1}$};

\node [node distance=2.5em, above of=f,rotate=-90] (b1) {$\prec$}; 
\node [node distance=1.5em, above of=b1,rotate=-90] (b2) {$...$};
\node [node distance=1.5em, above of=b2,rotate=-90] (b3) {$\prec$}; 
\node [draw,circle, fill=orange!20, node distance=2.5em, above of=b3] (b4) {$\varepsilon_1$};

\node [node distance=2.5em,rotate=25, above of=l, rotate=-90] (m) {$\prec$}; 
\node [node distance=2.5em,rotate=50, above of=l, rotate=-90] (n) {$\prec$}; 

\end{tikzpicture}
\captionsetup{justification=centering,margin=2cm}
\caption{Illustration of 1.|1.}
\label{phi_\ind1k_pic1}
\end{center}
\end{figure}

\begin{minipage}[t]{0.5\textwidth}
\begin{enumerate}
  \item[2.] If $\exists\ i$, s.t. $\varepsilon_i\nprec w_\ind\ \wedge\ \varepsilon_i\nsucc \vec{w}$, then define $$m=\max(i\ |1\leq i\leq k+1, \varepsilon_i\nprec w_\ind\ \wedge\ \varepsilon_i\nsucc \vec{w}),$$ then we have
\end{enumerate}
\begin{align*}&\cph \phi_{\ind,1^{k+1}}(\vec{w}\ |\vec{\varepsilon}\ )=\\&\cph =(\vec{w}\ ,\varepsilon_{m}; \varepsilon_1,...,\varepsilon_{m-1},\varepsilon_{m+1},..,\varepsilon_{k+1})\in  M^U_{\ind+1,1^{k}}.\end{align*}

\end{minipage}
\vrule \hspace{0.5cm} 
\begin{minipage}[t]{0.5\textwidth}
\begin{enumerate}
  \item[2.] For $(\vec{w}\ ,z ;\vec{\nu}\ )\in M^U_{\ind+1,1^{k}}$, we define $$j=\min(i\ |1\leq i\leq k, z\prec \nu_i),$$ then we have
\end{enumerate}
\begin{align*}&\cps\psi^2_{\ind,1^{k+1}}(\vec{w}\ ,z ;\vec{\nu}\ )=\\&\cps =(\vec{w}\ |\nu_1,...,\nu_{j-1},z,\nu_{j},..,\nu_k)\in P^U_\ind\times E_{k+1}^U.\end{align*}
\end{minipage}

\begin{figure}[H]
\begin{center}
\begin{tikzpicture} [>=latex,every node/.style={minimum width=3em, node distance=4em}, sh/.style={shade,shading=axis,left color=orange!20,right color=blue!20}]

\node [draw,circle, fill=blue!20] (b) {$\w_1$}; 
\node [node distance=5em,right of=b] (d) {$...$};
\node [draw,circle, fill=blue!20, node distance=5em,right of=d] (f1) {$\w_{\ind}$}; 
\node [node distance=2.5em, right of=f1] (k1) {$\preceq$}; 
\node [draw,circle, fill=blue!20!orange!20, node distance=2.5em,right of=k1] (f) [sh] {$\varepsilon_{m}|z$}; 
\node [node distance=2.5em, below of=f,rotate=-90] (k) {$\prec$}; 
\node [draw,circle, fill=orange!20, node distance=3em, below of=k] (l) {$\varepsilon_{m+1}|\nu_j$};
\node [node distance=2.75em, below of=l,rotate=-90] (zl) {$\prec$}; 
\node [node distance=1.25em, below of=zl,rotate=-90] (zl1) {$...$}; 
\node [node distance=1.25em, below of=zl1,rotate=-90] (zl2) {$\prec$}; 
\node [draw,circle, fill=orange!20, node distance=2.75em, below of=zl2] (zq1) {$\varepsilon_{k+1}|\nu_k$};

\node [node distance=2.5em, above of=f,rotate=-90] (b1) {$\prec$}; 
\node [node distance=1.5em, above of=b1,rotate=-90] (b2) {$...$};
\node [node distance=1.5em, above of=b2,rotate=-90] (b3) {$\prec$}; 
\node [draw,circle, fill=orange!20, node distance=2.5em, above of=b3] (b4) {$\varepsilon_1|\nu_1$};

\node [node distance=3em,rotate=30, above of=l, rotate=-90] (m) {$\prec$}; 
\node [node distance=3em,rotate=60, above of=l, rotate=-90] (n) {$\prec$}; 

\end{tikzpicture}
\captionsetup{justification=centering,margin=2cm}
\caption{Illustration of 2.|2.}
\label{phi_\ind1k_pic2}
\end{center}
\end{figure}
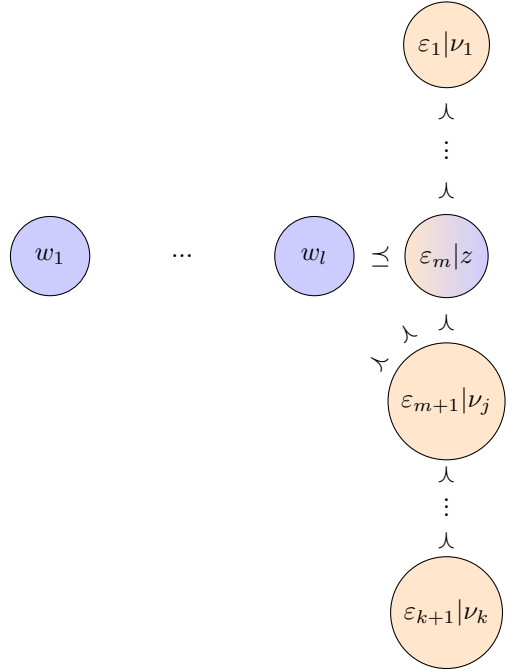
This completes the proof.
\end{proof}

\bigskip

\bigskip

Given a correct $\vs\in P^U_{\ind}$, let $$\theta= \theta(\vs)=\max(\{i<\ind|\; w_{i}\sim w_{i+1} \})\in\mathbb{N} \text{ and } J_{\ind-1}=(w_1,...,w_{\ind-1})\in P_{l-1}.$$ 

\begin{thm}\label{Thn2}
For natural $\ind\geq 2$, let $$M^U_{\ind,2} = \{(\vs\ \I q_0,q_1)\in P^U_{\ind.2}|\; J_{\ind-1}\prec q_0 \text{ and } w_{\ind}\prec q_1 \ \vee\ 
w_\theta\succ q_0 \text{ and }  w_{\theta+1}\succ q_1
\}.$$ Then,
$$m_{\ind,2}^U=\sum_{(\vs;q_0,q_1)\in M^U_{\ind,2}}w_1\cdot...\cdot w_{\ind}\cdot q_0\cdot q_1.$$

\end{thm}
\begin{remark}
According to Remark~\ref{c_m}, this implies $c_{n-2,2}(U)\geq 0$.
\end{remark}

  \begin{remark}
 There is a slightly more elegant version of
 $M^U_{\ind,2},$ which we will use in the future:
\[
   M^U_{\ind,2} = \{(\vs,q_0,q_1)\in P^U_{\ind,2}|\; (J_{\ind-1}\prec q_0 \wedge w_{\ind}\prec q_1) \ \vee\ 
(w_\theta\succ q_0 \wedge  w_{\theta+1}\succ q_1)
\}
\]

\end{remark}
Here is an illustration of an element of $M^U_{\ind,2}$:

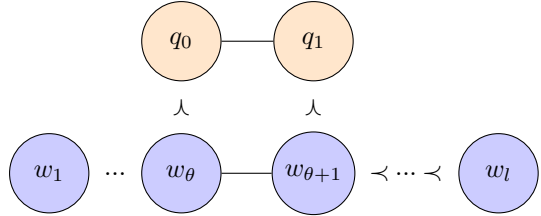
\begin{figure}[H]
\begin{center}
\begin{tikzpicture} [>=latex,every node/.style={minimum width=3em, node distance=4em}]

\node [draw,circle, fill=blue!20] (b) {$w_1$}; 
\node [node distance=2.5em,right of=b] (d) {$...$}; 
\node [draw,circle, fill=blue!20, node distance=2.5em,right of=d] (f) {$w_{\theta}$}; 
\node [draw,circle, fill=blue!20, node distance=5em,right of=f] (ff) {$w_{\theta+1}$}; 
\node [node distance=2.5em,right of=ff] (g) {$\prec$}; 
\node [node distance=1em,right of=g] (h) {$...$}; 
\node [node distance=1em,right of=h] (i) {$\prec$}; 
\node [draw,circle, fill=blue!20, node distance=2.5em,right of=i] (ii) {$w_{\ind}$}; 

\node [node distance=2.5em, above of=f,rotate=-90] (k) {$\prec$}; 
\node [draw,circle, fill=orange!20, node distance=2.5em, above of=k] (l) {$q_0$}; 
\node [node distance=2.5em, above of=ff,rotate=-90] (kk) {$\prec$}; 
\node [draw,circle, fill=orange!20, node distance=2.5em, above of=kk] (ll) {$q_1$}; 
\draw [-] (f) -- (ff);
\draw [-] (l) -- (ll);

\end{tikzpicture}
\caption{Illustration of  $M^U_{\ind,2}$.}
\label{Mn2_pic1}
\end{center}
\end{figure}

\begin{proof} We can write
\begin{equation}  \label{eqpn2}
   P^U_{\ind,2}\setminus M^U_{\ind,2} = \{(\vs\ \I q_0,q_1)\in P^U_{\ind.2}|\; (\vec{w}\prec q_0 \Rightarrow
(w_{\ind}\nprec q_1\wedge J_{\ind-1}\nprec q_1) \ \wedge\ 
( w_\theta\succ q_0\Rightarrow  w_{\theta+1}\not\succ q_1)
\}.
\end{equation}
The conditions in \eqref{eqpn2} have the form $A\wedge B$.
We begin with a few remarks.
\begin{enumerate}
\item Observe that the conditions of $A$ and $B$ are mutually exclusive, so we can
consider the two statements independently. 
\item Define $\tau=\max\{i\leq\ind|\; q_1\nprec w_{i}\vee w_i\sim w_{i+1}\}$. Note that it could happen that
$q_1\prec\theta$, but clearly $\tau\geq\theta$. 
\item For $\vec{u}\in P^U_{\ind+2}$, let $U_{\ind-1}=(u_1,...,u_{\ind-1})$ and $\breve{\theta}=\max(\{i<\ind|\; u_{i}\sim u_{i+1} \}).$
\end{enumerate}
To prove the theorem, we consider the following formula
$$m^U_{\ind,2}=p^U_\ind\cdot p^U_2-p^U_{\ind+2}.$$
and construct two injective maps. 

\bigskip

\begin{Parallel}[v]{0.48\textwidth}{0.48\textwidth}
\ParallelLText{
{\bf \Romannum{1}. We define $$\phi:P^U_{\ind,2}\setminus M^U_{\ind,2} \to P^U_{\ind+2}$$ as follows:}

Let $(\vs,q_0,q_1)\in P^U_{\ind,2}\setminus M^U_{\ind,2}.$

\begin{enumerate}
  \item[1.] If $\vec{w}\prec q_0$, then we define
$$\cph \phi(\vs,\vec{q}\ )=(... q_1,w_{\ind},q_0),\color{black}$$

which is in $P^U_{\ind+2},$ since $$w_{\ind}\nprec q_1\wedge J_{\ind-1}\nprec q_1.$$
\end{enumerate}
 }

\ParallelRText{
{\bf  \Romannum{2}. $$\psi:P^U_{\ind+2}\to P^U_{\ind,2}\setminus M^U_{\ind,2},$$ the inverse of map $\psi$:}

Let $(u_1,...,u_{\ind+2})\in P^U_{\ind+2}.$

\begin{enumerate}
  \item[1.] If $u_{\ind+2}\succ U_{\ind-1}$ and $u_{\ind+2}\succ u_{\ind+1}$, then 
$$\cps\psi(\vec{u}\ ) =(u_1,u_2,...,u_{\ind-1},u_{\ind+1}\I u_{\ind+2},u_{\ind})$$

\end{enumerate}

 }
\ParallelPar
\end{Parallel}

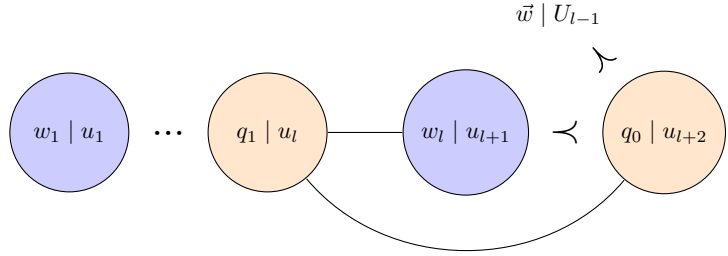
\begin{figure}[H]
\begin{center}
\begin{tikzpicture} [>=latex,every node/.style={minimum width=3em, node distance=4em},scale=1.5,transform shape]

\node [draw,circle, fill=blue!20] (b) {$\Scale[0.6] {w_1\mid u_1}$}; 
\node [node distance=2.5em,right of=b] (d) {$...$}; 
\node [draw,circle, fill=orange!20, node distance=2.5em,right of=d] (f) {$\Scale[0.6] {q_1\mid u_{\ind}}$}; 
\node [draw,circle, fill=blue!20, node distance=5em,right of=f] (ii) {$\Scale[0.6] {w_{\ind}\mid u_{\ind+1}}$}; 
\node [node distance=2.5em,right of=ii] (iii) {$\prec$}; 
\node [draw,circle, fill=orange!20, node distance=2.5em,right of=iii] (j) {$\Scale[0.6] {q_0\mid u_{\ind+2}}$}; 
\node [node distance=2.5em, rotate=40 ,above of=j,rotate=-90] (k) {$\prec$}; 
\node [node distance=1.5em, above left of=k] (n) {$\Scale[0.6] {\vec{w}\mid U_{\ind-1}}$}; 
\draw [-] (f) -- (ii);
\draw [-] (f) to [out=-50,in=-130] (j);

\end{tikzpicture}
\captionsetup{justification=centering,margin=2cm}
\caption{Illustration of $\phi_{\ind|2}$, if $q_0\succ \vec{w}$.\\Illustration of $\psi_{\ind|2}$, if $u_{\ind+2}\succ U_{\ind-1}$ and $u_{\ind+2}\succ u_{\ind+1}$.}
\label{phi_\ind2_pi1}
\end{center}
\end{figure}


\begin{Parallel}[v]{0.48\textwidth}{0.48\textwidth}
\ParallelLText{
\
\begin{enumerate}
 \item[2.]  $q_0\prec w_\tau$ implies $A$, and implies $\tau>\theta$; then  we define
$$\cph\phi(\vs,\vec{q}\ )=(... w_{\tau},q_1,q_0,w_{\tau+1},...).$$

Here and below, $\tau+1$ might be equal to $\ind+1$, in which case we
simply omit $w_{\tau+1}$.
\end{enumerate}
 }
\ParallelRText{
\ 
\begin{enumerate}

\item[2.] if $u_{\breve{\theta}+1}\prec u_{\breve{\theta}-1}$, then 

$$\cps\psi(\vec{u}\ ) =(u_1,u_2,...,u_{\breve{\theta}-1},u_{\breve{\theta}+2},...,u_{\ind+2}\I u_{\breve{\theta}+1},u_{\breve{\theta}}),$$

Here and below, $\breve{\theta}+1$ might be equal to $\ind+2$, in which case we
simply omit $u_{\breve{\theta}+2}$.

\end{enumerate}
 }
\ParallelPar
\end{Parallel}

\begin{figure}[H]
\begin{center}
\begin{tikzpicture} [>=latex,every node/.style={minimum width=3em, node distance=4em},scale=1.5,transform shape]

\node [draw,circle, fill=blue!20] (b) {$\Scale[0.6]{w_1\mid u_1}$}; 
\node [node distance=2.5em,right of=b] (d) {$...$}; 

\node [draw,circle, fill=blue!20, node distance=2.5em,right of=d] (dd) {$\Scale[0.6]{w_\tau\mid u_{\breve{\theta}-1}}$}; 
\node [node distance=2.5em, rotate=-40 ,above of=dd,rotate=-90] (k) {$\prec$}; 
\node [node distance=1.5em, above right of=k] (n) {$\Scale[0.6]{q_0\mid u_{\breve{\theta}+1}}$};
\node [draw,circle, fill=orange!20, node distance=5em,right of=dd] (f) {$\Scale[0.6]{q_1\mid u_{\breve{\theta}}}$}; 
\node [draw,circle, fill=orange!20, node distance=5em,right of=f] (fff) {$\Scale[0.6]{q_0\mid u_{\breve{\theta}+1}}$};  
\node [node distance=2.5em,right of=fff] (g) {$\prec$}; 
\node [draw,circle, fill=blue!20, node distance=2.5em,right of=g] (t1) {$\Scale[0.6]{w_{\tau+1}\mid u_{\breve{\theta}+2}}$};
\node [node distance=2.5em,right of=t1] (t2) {$\prec$}; 
\node [node distance=1em,right of=t2] (h) {$...$}; 

\draw [-] (f) -- (fff);
\draw [-] (dd) -- (f);
\end{tikzpicture}
\captionsetup{justification=centering,margin=2cm}
\caption{Illustration of $\phi_{\ind|2}$, if $q_0\prec w_\tau$.\\Illustration of $\psi_{\ind|2}$, if $u_{\breve{\theta}+1}\prec u_{\breve{\theta}-1}$.}
\label{phi_\ind2_pi2}
\end{center}
\end{figure}

\begin{Parallel}[v]{0.48\textwidth}{0.48\textwidth}

\ParallelLText{
\ 
\begin{enumerate}

  \item[3.]  Finally, if $\vec{w}\nprec q_0$ and $q_0\nprec w_\tau$, then
$$\cph\phi(\vs,\vec{q}\ )=(... w_{\tau},q_0,q_1,w_{\tau+1},...).\color{black}$$

Here and below, $\tau+1$ might be equal to $\ind+1$, in which case we
simply omit $w_{\tau+1}$.
\end{enumerate}
 }

\ParallelRText{
\ 
\begin{enumerate}

\item[3.] if $u_{\breve{\theta}+1}\nprec u_{\breve{\theta}-1}$, then $$\cps\psi(\vec{u}\ ) =(u_1,u_2,...,u_{\breve{\theta}-1},u_{\breve{\theta}+2},...,u_{\ind+2}\I u_{\breve{\theta}},u_{\breve{\theta}+1}).$$

Here and below, $\breve{\theta}+1$ might be equal to $\ind+2$, in which case we
simply omit $u_{\breve{\theta}+2}$.

\end{enumerate}

 }
\ParallelPar
\end{Parallel}

\begin{figure}[H]
\begin{center}
\begin{tikzpicture} [>=latex,every node/.style={minimum width=3em, node distance=4em},scale=1.5,transform shape]

\node [draw,circle, fill=blue!20] (b) {$\Scale[0.6]{w_1\mid u_1}$}; 
\node [node distance=2.5em,right of=b] (d) {$...$}; 

\node [draw,circle, fill=blue!20, node distance=2.5em,right of=d] (dd) {$\Scale[0.6]{w_\tau\mid u_{\breve{\theta}-1}}$}; 
\node [node distance=2.5em, rotate=-40 ,above of=dd,rotate=-90] (k) {$\nprec$}; 
\node [node distance=1.5em, above right of=k] (n) {$\Scale[0.6]{q_1\mid u_{\breve{\theta}+1}}$};
\node [node distance=2.5em,right of=dd] (dd1) {$\preceq$}; 
\node [draw,circle, fill=orange!20, node distance=2.5em,right of=dd1] (f) {$\Scale[0.6]{q_0\mid u_{\breve{\theta}}}$}; 
\node [draw,circle, fill=orange!20, node distance=5em,right of=f] (fff) {$\Scale[0.6]{q_1\mid u_{\breve{\theta}+1}}$}; 
\node [node distance=2.5em,right of=fff] (g) {$\prec$}; 
\node [draw,circle, fill=blue!20, node distance=2.5em,right of=g] (t1) {$\Scale[0.6]{w_{\tau+1}\mid u_{\breve{\theta}+2}}$};
\node [node distance=2.5em,right of=t1] (t2) {$\prec$}; 
\node [node distance=1em,right of=t2] (h) {$...$};

\draw [-] (f) -- (fff);
\end{tikzpicture}
\captionsetup{justification=centering,margin=2cm}
\caption{Illustration of $\phi_{\ind|2}$, if $q_0\nprec w_\tau$.\\Illustration of $\psi_{\ind|2}$, if $u_{\breve{\theta}+1}\nprec u_{\breve{\theta}-1}$.}
\label{psi_\ind2_pi3}
\end{center}
\end{figure}

\end{proof}

To find a combinatorial interpretation of $m^U_{\ind,2,1}$, we construct a bijection between the right and left hand sides of the following equality:

\begin{equation} \label{n21}
p^U_{\ind}*m^U_{2,1}=m^U_{\ind+2,1}+m^U_{\ind+1,2}+m^U_{\ind,2,1}.
\end{equation}

This formula and its proof are similar to the previous one.

\begin{thm}\label{thmn21}
For natural $\ind$, let
\begin{equation}
\begin{aligned}
M^U_{\ind,2,1}= &\{ (\vec{w}\ |\vec{q}\ |z)\in P_{\ind,2,1}|\; (\vec{w}\ ;\vec{q}\ )\in M^U_{\ind,2},\ (\vec{w}\ ;z)\in M^U_{\ind,1},\ (\vec{q}\ ;z)\in M^U_{2,1},  \}\\ &  \cup \{ (\vec{w}\ |\vec{q}\ |z)\in P_{\ind,2,1}|\; s_\ind\succ z\succ\vec{q},\ \exists \gamma \in \mathbb{N}, \text{ such that } \theta < \gamma<\ind,\ s_{ \gamma}\sim z, \ s_{ \gamma}\succ q_2 \}\\ &  \cup \{ (\vec{w}\ |\vec{q}\ |z)\in P_{\ind,2,1}|\; s_\ind\succ z\succ q_2,\ \exists \gamma \in \mathbb{N}, \text{ such that } \theta < \gamma<\ind,\ s_{\gamma}\sim z,\ z\sim q_1,\ s_{ \gamma}\succ q_1 \},
\end{aligned}
\end{equation}
Then,
$$m_{\ind,2,1}^U=\sum_{(\vs;q_0,q_1;z)\in M^U_{\ind,2,1}}w_1\cdot...\cdot w_{\ind}\cdot q_0\cdot q_1\cdot z.$$
\end{thm}
\begin{remark}
According to Remark~\ref{c_m}, this implies $c_{n-3,2,1}(U)\geq 0$.
\end{remark}

Let us explain the meaning of $M_{\ind,2,1}^U$. This theorem states that in addition to combinations of pairwise comparable corrects of lengths $\ind$, 2 and 1 we have two more cases:

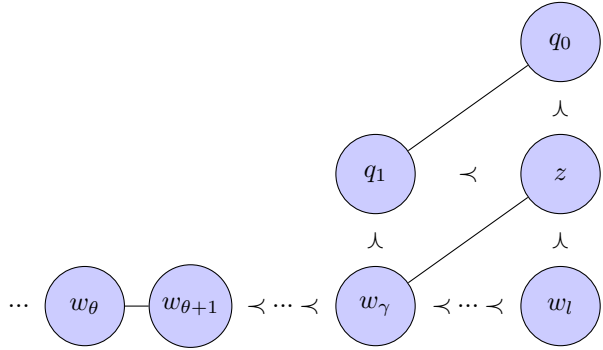
\begin{figure}[H]
\begin{center}
\begin{tikzpicture} [>=latex,every node/.style={minimum width=3em, node distance=4em}]

\node  (start) {$...$}; 
\node [draw,circle, fill=blue!20, node distance=2.5em, right of=start] (a)  {$w_{\theta}$};
\node [draw,circle, fill=blue!20, right of=a] (b) {$w_{\theta+1}$}; 
\node [node distance=2.5em,right of=b] (c) {$\prec$}; 
\node [node distance=1em,right of=c] (d) {$...$}; 
\node [node distance=1em,right of=d] (e) {$\prec$}; 
\node [draw,circle, fill=blue!20, node distance=2.5em,right of=e] (f) {$w_{\gamma}$}; 
\node [node distance=2.5em,right of=f] (g) {$\prec$}; 
\node [node distance=1em,right of=g] (h) {$...$}; 
\node [node distance=1em,right of=h] (i) {$\prec$}; 
\node [draw,circle, fill=blue!20, node distance=2.5em,right of=i] (j) {$w_{\ind}$}; 
\node [node distance=2.5em, above of=j,rotate=-90] (k) {$\prec$}; 
\node [draw,circle, fill=blue!20, node distance=2.5em,above of=k] (l) {$z$}; 
\node [node distance=3.5em, left of=l] (m) {$\prec$}; 
\node [draw,circle, fill=blue!20, node distance=3.5em,left of=m] (n) {$q_1$}; 
\node [node distance=2.5em, above of=f,rotate=-90] (o) {$\prec$}; 
\node [node distance=2.5em, above of=l,rotate=-90] (p) {$\prec$}; 
\node [draw,circle, fill=blue!20, node distance=2.5em,above of=p] (q) {$q_0$}; 
\draw [-] (b) -- (a);
\draw [-] (l) -- (f);
\draw [-] (n) -- (q);

\end{tikzpicture}
\captionsetup{justification=centering,margin=2cm}
\caption{First exceptional element of $M^U_{\ind,2,1}$}
\label{Firstn21}
\end{center}

\end{figure}

\begin{figure}[H]
\begin{center}
\begin{tikzpicture} [>=latex,every node/.style={minimum width=3em, node distance=4em}]

\node  (start) {$...$}; 
\node [draw,circle, fill=blue!20, node distance=2.5em, right of=start] (a)  {$w_{\theta}$};
\node [draw,circle, fill=blue!20, right of=a] (b) {$w_{\theta+1}$}; 
\node [node distance=2.5em,right of=b] (c) {$\prec$}; 
\node [node distance=1em,right of=c] (d) {$...$}; 
\node [node distance=1em,right of=d] (e) {$\prec$}; 
\node [draw,circle, fill=blue!20, node distance=2.5em,right of=e] (f) {$w_{\gamma}$}; 
\node [node distance=2.5em,right of=f] (g) {$\prec$}; 
\node [node distance=1em,right of=g] (h) {$...$}; 
\node [node distance=1em,right of=h] (i) {$\prec$}; 
\node [draw,circle, fill=blue!20, node distance=2.5em,right of=i] (j) {$w_{\ind}$}; 
\node [node distance=2.5em, above of=j,rotate=-90] (k) {$\prec$}; 
\node [draw,circle, fill=blue!20, node distance=2.5em, above of=k] (l) {$q_1$}; 

\node [draw,circle, fill=blue!20, node distance=7em,left of=l] (n) {$q_0$}; 
\node [node distance=2.5em, above of=f,rotate=-90] (o) {$\prec$}; 
\node [node distance=2.5em, above of=l,rotate=-90] (p) {$\prec$}; 
\node [draw,circle, fill=blue!20, node distance=2.5em,above of=p] (q) {$z$}; 
\draw [-] (b) -- (a);
\draw [-] (l) -- (f);
\draw [-] (l) -- (n);
\draw [-] (n) -- (q);

\end{tikzpicture}
\captionsetup{justification=centering,margin=2cm}
\caption{Second exceptional element of $M^U_{\ind,2,1}$}
\label{Secondn21}
\end{center}
\end{figure}
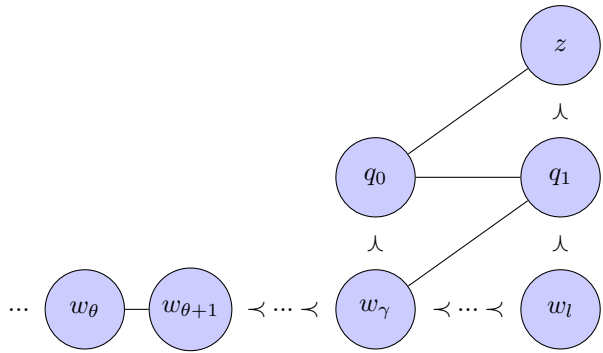

\begin{proof}

To prove Theorem \ref{thmn21} using Formula \ref{n21}, we construct  the maps $\varphi_{\ind|2,1}$ and $\psi_{\ind|2,1}$.

\bigskip

\begin{minipage}[t]{0.5\textwidth}

{\bf \Romannum{1}}. We construct the map from the left hand side to the right hand side  $$\phi_{\ind|2,1}: P^U_{\ind}\times M^U_{2,1} \to M_{\ind+2,1}^U\sqcup M_{\ind+1,2}^U\sqcup M^U_{\ind,2,1}.$$

\smallskip Let us take $$(\vec{\w}\ |q_1,q_2\ ;z)\in P^U_{\ind}\times M^U_{2,1}.$$ 

Let $$\theta=\max(\{i<\ind|\; \w_{i}\sim \w_{i+1} \}).$$

We will use this $\theta$ for $\w$ on the right hand side as well.

\begin{enumerate}[start=1,leftmargin=10pt,label*={\arabic*}]
\item[1.] If $z\succ\vec{q}$.
 \begin{enumerate}[start=1,leftmargin=5pt,label*={\arabic*}]
   \item[11.]  If $z\succ\vec{\w}$.
    \begin{enumerate}[start=1,leftmargin=5pt,label*={\arabic*}]
      \item[111.]  If $(\vec{\w}\ ;\vec{q}\ )\in M^U_{\ind,2}$, then 
  $$\cph \phi_{\ind|2,1}(\vec{\w}\ |\vec{q}\ ;z)=(\vec{\w}\ ;\vec{q}\ ;z)\in M^U_{\ind,2,1}.$$

\end{enumerate}
\end{enumerate}
\end{enumerate}

\end{minipage}
\vrule \hspace{0.5cm} 
\begin{minipage}[t]{0.5\textwidth}

{\bf  \Romannum{2}.} We construct $$\psi^1_{\ind|2,1}: M^U_{\ind+2,1} \to P^U_\ind\times M^U_{2,1},$$    $$\psi^2_{\ind|2,1}: M^U_{\ind+1,2} \to P^U_\ind\times M^U_{2,1},$$   $$\psi^3_{\ind|2,1}: M^U_{\ind,2,1} \to P^U_\ind\times M^U_{2,1}.$$ 
We take $$(\vec{\wii}\ ;\xix) \in M^U_{\ind+2,1};$$ $$(\vec{\w},\xix\ ;\ \kap_0,\kap_1) \in M^U_{\ind+1,2};$$ $$(\vec{\w}\ ;\kap_0,\kap_1;z) \in M^U_{\ind,2,1}.$$


Let $\breve{\theta}=\max(\{i<\ind+2|\ \wii_i\sim \wii_{i+1}\})$.

\begin{enumerate}[start=1,leftmargin=5pt,label*={\arabic*}]
\item [1.] If $(\vec{\w}\ ;\kap_0,\kap_1; z) \in M^U_{\ind,2,1}$ and $z\succ\vec{\kap}$ and $z\succ\vec{\w}$, then
 $$\cps \psi^3_{\ind|2,1}(\vec{\w}\ ;\kap_0,\kap_1;z)=(\vec{\w}\ \I\ \kap_0,\kap_1\ ; z).$$ 

\end{enumerate}
\end{minipage}

\begin{center} In this case we have 2 illustrations:\end{center}
\bigskip

\begin{figure}[H]
\begin{center}
\begin{tikzpicture} [>=latex,every node/.style={minimum width=3em, node distance=4em}]

\node [draw,circle, fill=blue!20] (b) {$\w_1$}; 
\node [node distance=2.5em,right of=b] (d) {$...$}; 
\node [draw,circle, fill=blue!20, node distance=2.5em,right of=d] (f) {$\w_{\ind-1}$}; 
\node [node distance=2.5em,right of=f] (i) {$\preceq$}; 
\node [draw,circle, fill=blue!20, node distance=2.5em,right of=i] (j) {$\w_\ind$};
\node [node distance=2.5em, below of=f,rotate=-90] (k) {$\prec$}; 
\node [draw,circle, fill=orange!20, node distance=2.5em, below of=k] (l) {$q_0$}; 
\node [node distance=2.5em,rotate=25, above of=l, rotate=-90] (m) {$\prec$}; 
\node [node distance=2.5em,rotate=50, above of=l, rotate=-90] (n) {$\prec$}; 
\node [node distance=2.5em, below of=j,rotate=-90] (o) {$\prec$}; 
\node [draw,circle, fill=orange!20, node distance=2.5em, below of=o] (p) {$q_1$};
\node [node distance=2.5em, below of=p,rotate=-90] (oo) {$\prec$}; 
\node [draw,circle, fill=orange!50, node distance=2.5em, below of=oo] (pp) {$z$};
\node [node distance=2.5em,rotate=45, above of=pp, rotate=-90] (nn) {$\prec$}; 
\draw [->] (l) -- (p);

\end{tikzpicture}
\captionsetup{justification=centering,margin=2cm}
\caption{Illustration of $\phi_{\ind|2,1}$, case 111., and $\psi^3_{\ind|2,1}$, case 1.}
\label{phi_\ind21_pic0a}
\end{center}
\end{figure}

\begin{figure}[H]
\begin{center}
\begin{tikzpicture} [>=latex,every node/.style={minimum width=3em, node distance=4em}]

\node [draw,circle, fill=blue!20] (b) {$\w_1$}; 
\node [node distance=2.5em,right of=b] (d) {$...$}; 
\node [draw,circle, fill=blue!20, node distance=2.5em,right of=d] (f) {$\w_{\theta}$}; 
\node [draw,circle, fill=blue!20, node distance=5em,right of=f] (ff) {$\w_{\theta+1}$}; 
\node [node distance=2.5em,right of=ff] (g) {$\prec$}; 
\node [node distance=1em,right of=g] (h) {$...$}; 
\node [node distance=1em,right of=h] (i) {$\prec$}; 
\node [draw,circle, fill=blue!20, node distance=2.5em,right of=i] (ii) {$\w_{\ind}$}; 
\node [node distance=2.5em, below of=ii,rotate=-90] (iii) {$\prec$}; 
\node [draw,circle, fill=orange!50, node distance=2.5em,below of=iii] (j) {$z$}; 
\node [node distance=2.5em,rotate=25, above of=j, rotate=-90] (m) {$\prec$}; 
\node [node distance=2.5em,rotate=50, above of=j, rotate=-90] (n) {$\prec$};

\node [node distance=2.5em, above of=f,rotate=-90] (k) {$\prec$}; 
\node [draw,circle, fill=orange!20, node distance=2.5em, above of=k] (l) {$q_0$}; 
\node [node distance=2.5em, above of=ff,rotate=-90] (kk) {$\prec$}; 
\node [draw,circle, fill=orange!20, node distance=2.5em, above of=kk] (ll) {$q_1$}; 
\draw [-] (f) -- (ff);
\draw [->] (l) -- (ll);

\end{tikzpicture}
\captionsetup{justification=centering,margin=2cm}
\caption{Illustration of $\phi_{\ind|2,1}$, case 111., and $\psi^3_{\ind|2,1}$, case 1.}
\label{phi_\ind21_pic0b}
\end{center}
\end{figure}

\begin{minipage}[t]{0.5\textwidth}

    \begin{enumerate}[start=112,leftmargin=15pt,label*={\arabic*}]
      \item[112.] If $(\vec{\w}\ |\vec{q}\ )\notin M^U_{\ind,2}$, then using $$\phi_{\ind|2}(\vec{\w}\ |\vec{q}\ )\in M^U_{\ind+2}$$ we have
$$\cph \phi_{\ind|2,1}(\vec{\w}\ |\vec{q}\ |z)=(\phi_{\ind|2}(\vec{\w}\ |\vec{q}\ )\ ;z)\in M^U_{\ind+2,1}$$
    \end{enumerate}

\end{minipage}
\vrule \hspace{0.5cm} 
\begin{minipage}[t]{0.5\textwidth}

\begin{enumerate}[start=2,leftmargin=5pt,label*={\arabic*}]
\item[2.] If $(\vec{\wii}\ ;\xix) \in M^U_{\ind+2,1}$ and $\xix\succ\vec{\wii}$, then using $$\psi_{\ind|2}(\vec{\wii}\ )\in M^U_{\ind,2}$$ we have
$$\cps \psi^1_{\ind|2,1}(\vec{\wii}\ ;z)=(\psi_{\ind|2}(\vec{\wii}\ );\ \xix).$$

\end{enumerate}

\end{minipage}

\begin{center} Here, we provide illustrations for the right hand side, see Theorem \ref{Thn2}, where $\phi_{\ind|2}$ and $\psi_{\ind|2}$ are defined, for more details. \end{center}

\bigskip

\begin{figure}[H]
\begin{center}
\begin{tikzpicture} [>=latex,every node/.style={minimum width=3em, node distance=4em}]

\node [draw,circle, fill=blue!20] (start)  {$\wii_1$};
\node  [node distance=2.5em, right of=start] (start1) {$...$}; 
\node [draw,circle, fill=orange!20, node distance=2.5em, right of=start1] (a)  {$\wii_{\theta}$};
\node [draw,circle, fill=orange!20, right of=a] (b) {$\wii_{\theta+1}$}; 
\node [node distance=2.5em,right of=b] (c) {$\prec$}; 
\node [draw,circle, fill=blue!20, node distance=2.5em, right of=c] (cc) {$\wii_{\theta+2}$}; 
\node [node distance=2.5em,right of=cc] (ccc) {$\prec$}; 
\node [node distance=1em,right of=ccc] (d) {$...$}; 
\node [node distance=1em,right of=d] (e) {$\prec$}; 

\node [draw,circle, fill=blue!20, node distance=2.5em,right of=e] (fff) {$\wii_{\ind+1}$}; 

\node [node distance=2.5em,right of=fff] (i) {$\prec$}; 
\node [draw,circle, fill=blue!20, node distance=2.5em,right of=i] (j) {$\wii_{\ind+2}$}; 
\node [node distance=2.5em, below of=j,rotate=-90] (k) {$\prec$}; 
\node [draw,circle, fill=orange!50, node distance=2.5em, below of=k] (l) {$\xix$}; 
\node [node distance=2.5em,rotate=25, above of=l, rotate=-90] (m) {$\prec$}; 
\node [node distance=2.5em,rotate=50, above of=l, rotate=-90] (m) {$\prec$}; 
\node [node distance=2.5em,rotate=75, above of=l, rotate=-90] (m) {$\prec$}; 
\draw [<->] (b) -- (a);

\end{tikzpicture}
\captionsetup{justification=centering,margin=2cm}
\caption{Illustration of $\psi^1_{\ind|2,1}$, when $\xix\succ\vec{\wii}$, general case }
\label{psi1_\ind21_pic1}
\end{center}
\end{figure}
Note that the picture above illustrates most of the cases, except the special one, when $\wii_{\ind+2}~\succ~\wii_{\ind+1}$, $\wii_{\ind+2}\succ J_{\ind-1}$ and $\wii_{\ind+2}\sim \wii_{\ind}$. In this case map $\psi_{\ind|2}$ takes out $\wii_{\ind+2}$ and $\wii_{\ind}$:

\begin{figure}[H]
\begin{center}
\begin{tikzpicture} [>=latex,every node/.style={minimum width=3em, node distance=4em}]

\node [draw,circle, fill=blue!20] (b) {$\wii_1$}; 
\node [node distance=2.5em,right of=b] (d) {$...$}; 
\node [draw,circle, fill=orange!20, node distance=2.5em,right of=d] (f) {$\wii_\ind$}; 
\node [draw,circle, fill=blue!20, node distance=5em,right of=f] (ii) {$\wii_{\ind+1}$}; 
\node [node distance=2.5em,right of=ii] (iii) {$\prec$}; 
\node [draw,circle, fill=orange!20, node distance=2.5em,right of=iii] (j) {$\wii_{\ind+2}$}; 
\node [node distance=2.5em, below of=j,rotate=-90] (k) {$\prec$}; 
\node [draw,circle, fill=orange!50, node distance=2.5em, below of=k] (l) {$\xix$}; 
\node [node distance=2.5em,rotate=25, above of=l, rotate=-90] (m) {$\prec$}; 
\node [node distance=2.5em,rotate=50, above of=l, rotate=-90] (m) {$\prec$}; 
\node [node distance=2.5em,rotate=75, above of=l, rotate=-90] (m) {$\prec$}; 
\draw [-] (f) -- (ii);
\draw [->] (j) to [out=130,in=50] (f);

\end{tikzpicture}
\captionsetup{justification=centering,margin=2cm}
\caption{Illustration of $\psi^1_{\ind|2,1}$ (special case),\\ if $\wii_{\ind+2}~\succ~\wii_{\ind+1}$, $\wii_{\ind+2}\succ \{\wii_{1},...,\wii_{\ind-1}\}$ and $\xix\succ\vec{\wii}$.}
\label{psi1_\ind21_pic2}
\end{center}
\end{figure}

\begin{minipage}[t]{0.5\textwidth}

\begin{enumerate}[start=12,leftmargin=10pt,label*={\arabic*}]
   \item [12.] If $z\nprec \w_\ind$ and $z\nsucc\vec{w}$, then denote 
 \begin{equation*}
    (\hat{w}_\theta,\hat{w}_{\theta+1})=
    \begin{cases}
      (w_\ind,z), & \text{if}\ w_\ind\sim z, \\
      (\w_\theta,\w_{\theta+1}), & z\succ w_\ind.
    \end{cases}
  \end{equation*}
     \begin{enumerate}[leftmargin=5pt,label*={\arabic*}]
        \item [121.] If $\hat{w}_\theta\succ q_0$ and $\hat{w}_{\theta+1}\succ q_1$, then  
$$\cph\phi_{\ind|2,1}(\vec{\w}\ |q_0,q_1;z)=(\vec{\w},z\ ;q_0,q_1)\in M^U_{\ind+1,2}.$$
\end{enumerate}
\end{enumerate}

\end{minipage}
\vrule \hspace{0.5cm} 
\begin{minipage}[t]{0.5\textwidth}

\begin{enumerate}[start=3,leftmargin=5pt,label*={\arabic*}]
\item [3.]  For $(\vec{\w},\xix\ ;\kap_0,\kap_1) \in M^U_{\ind+1,2}$, such that

 ($\xix\succ w_\ind,$ $\w_{\theta}\succ \kap_0,$ $\w_{\theta+1}\succ\kap_1,$ and $\xix\succ \kap_0$), \\or

($\xix\sim w_\ind,$ $\w_{\ind}\succ \kap_0,$ $\xix\succ\kap_1,$ and $\xix\succ \kap_0$), \\we have:
$$\cps \psi^2_{\ind|2,1}(\vec{\w},\xix\ ;\kap_0,\kap_1)=(\vec{\w}\ |\kap_0,\kap_1;\xix).$$

\end{enumerate}

\end{minipage}
\bigskip

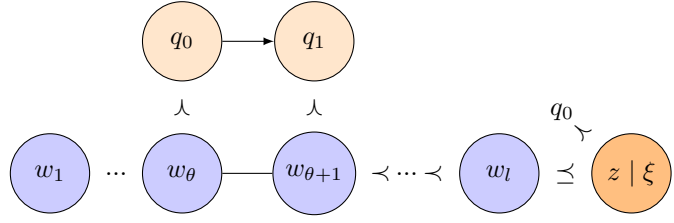
\begin{figure}[H]
\begin{center}
\begin{tikzpicture} [>=latex,every node/.style={minimum width=3em, node distance=4em}]

\node [draw,circle, fill=blue!20] (b) {$\w_1$}; 
\node [node distance=2.5em,right of=b] (d) {$...$}; 
\node [draw,circle, fill=blue!20, node distance=2.5em,right of=d] (f) {$\w_{\theta}$}; 
\node [draw,circle, fill=blue!20, node distance=5em,right of=f] (ff) {$\w_{\theta+1}$}; 
\node [node distance=2.5em,right of=ff] (g) {$\prec$}; 
\node [node distance=1em,right of=g] (h) {$...$}; 
\node [node distance=1em,right of=h] (i) {$\prec$}; 
\node [draw,circle, fill=blue!20, node distance=2.5em,right of=i] (ii) {$\w_{\ind}$}; 
\node [node distance=2.5em,right of=ii] (iii) {$\preceq$}; 
\node [draw,circle, fill=orange!50, node distance=2.5em,right of=iii] (j) {$z\mid\xix$}; 

\node [node distance=2.5em,rotate=50, above of=j, rotate=-90] (ar2) {$\prec$}; 
\node [node distance=1em, above left of=ar2] (arr2) {$q_0$}; 

\node [node distance=2.5em, above of=f,rotate=-90] (k) {$\prec$}; 
\node [draw,circle, fill=orange!20, node distance=2.5em, above of=k] (l) {$q_0$}; 
\node [node distance=2.5em, above of=ff,rotate=-90] (kk) {$\prec$}; 
\node [draw,circle, fill=orange!20, node distance=2.5em, above of=kk] (ll) {$q_1$}; 
\draw [-] (f) -- (ff);
\draw [->] (l) -- (ll);

\end{tikzpicture}
\captionsetup{justification=centering,margin=2cm}
\caption{ Illustration of 121.| 3.}
\label{phi_\ind21_pic1}
\end{center}
\end{figure}

\bigskip

\begin{minipage}[t]{0.5\textwidth}

     \begin{enumerate}[start=122,leftmargin=15pt,label*={\arabic*}]
        \item[122.] If $\hat{w}_{\theta+1}\nsucc q_1$ or $\hat{w}_{\theta}\nsucc q_0$, we take $$\tau =\max(\{i<\ind| \w_{i}\nsucc q_1 \vee \w_i\sim \w_{i+1}\vee \w_i\sim z\})$$ (note that $\tau\geq\theta$), and insert $q_0$ or $q_1$   after it: 
                \begin{enumerate}[start=1,leftmargin=5pt,label*={\arabic*}]
                   \item [1221.] If $q_0\prec \w_\tau$, then 
\end{enumerate}
\begin{multline*}\cph \phi_{\ind|2,1}(\vec{\w}\ |\vec{q}\ ;z)=\\ \cph= (\w_1,...,\w_\tau,q_1,\w_{\tau+1},...,\w_\ind,z\ ;q_0)\in M^U_{\ind+2,1} \end{multline*}

\end{enumerate}

\end{minipage}
\vrule \hspace{0.5cm} 
\begin{minipage}[t]{0.5\textwidth}

\begin{enumerate}[start=4,leftmargin=5pt,label*={\arabic*}]
\item For $(\vec{\wii}\ ;\xix)\in M^U_{\ind+2,1},$ such that $\xix\prec\wii_{\ind+2}$, define  $$\eta=\max(0,\{i\geq\breve{\theta}| \wii_i\sim \xix\}).$$
          \begin{enumerate} [leftmargin=5pt,label*={\arabic*}]
              \item[41.] If $\eta>0$, then we take out  $\wii_{\eta},\ \wii_{\ind+2}$ and $\xix$:
\begin{enumerate}  [leftmargin=5pt,label*={\arabic*}]
  \item [411.] If $\eta=\breve{\theta}+1$ and $\wii_{\breve{\theta}}\succ \xix$, then
\end{enumerate}
\setstretch{1}
\begin{multline*}\cps\psi^1_{\ind|2,1}(\vec{\wii}\ ;\xix) = (\wii_1,...,\wii_{\eta-1},\wii_{\eta+1},...,\wii_{\ind+1}\ |\xix,\wii_{\eta};\wii_{\ind+2}).\end{multline*}

\end{enumerate}
\end{enumerate}
\end{minipage}
\bigskip

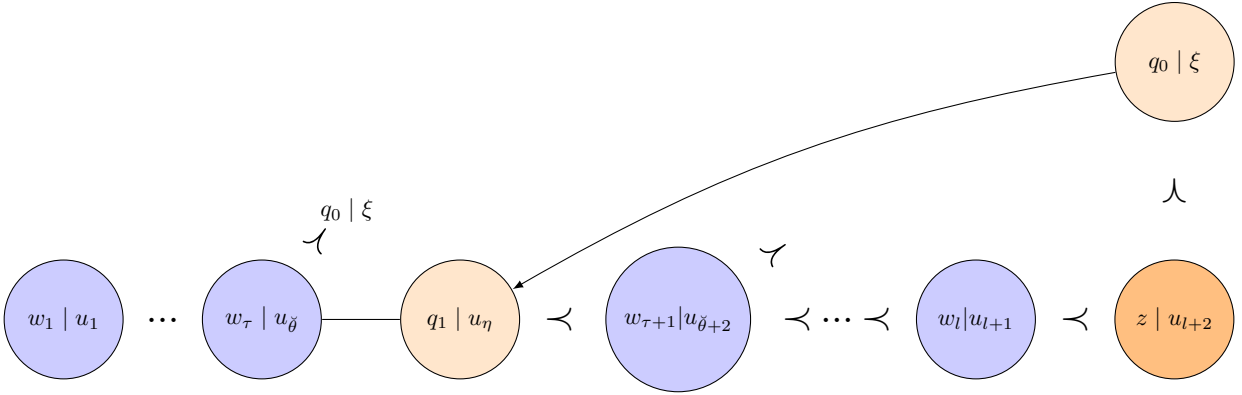
\begin{figure}[H]
\begin{center}
\begin{tikzpicture} [>=latex,every node/.style={minimum width=3em, node distance=4em},scale=1.5,transform shape]

\node [draw,circle, fill=blue!20] (b) {$\Scale[0.6]{\w_1\mid\wii_1}$}; 
\node [node distance=2.5em,right of=b] (d) {$...$}; 
\node [draw,circle, fill=blue!20, node distance=2.5em,right of=d] (e) {$\Scale[0.6]{\w_{\tau}\mid\wii_{\breve{\theta}}}$};
\node [node distance=2.5em,rotate=-35, above of=e, rotate=-90] (ee) {$\prec$};
\node [node distance=1em, above right of=ee] (eee) {$\Scale[0.6]{q_0\mid\xix}$};

\node [draw,circle, fill=orange!20, node distance=5em,right of=e] (f) {$\Scale[0.6]{q_1\mid\wii_{\eta}}$}; 
\node [node distance=2.5em,right of=f] (ff) {$\prec$}; 
\node [draw,circle, fill=blue!20, node distance=5.5em,right of=f] (fff) {$\Scale[0.6]{\w_{\tau+1}|\wii_{\breve{\theta}+2}}$}; 
\node [node distance=3em,rotate=-55, above of=fff, rotate=-90] (ffff) {$\prec$}; 
\node [node distance=3em,right of=fff] (g) {$\prec$}; 
\node [node distance=1em,right of=g] (h) {$...$}; 
\node [node distance=1em,right of=h] (i) {$\prec$}; 
\node [draw,circle, fill=blue!20, node distance=2.5em,right of=i] (ii) {$\Scale[0.6]{\w_\ind|\wii_{\ind+1}}$}; 
\node [node distance=2.5em,right of=ii] (iii) {$\prec$}; 
\node [draw,circle, fill=orange!50, node distance=2.5em,right of=iii] (j) {$\Scale[0.6]{ z\mid\wii_{\ind+2}}$}; 
\node [node distance=3.25em, above of=j,rotate=-90] (k) {$\prec$}; 
\node [draw,circle, fill=orange!20, node distance=3.25em, above of=k] (l) {$\Scale[0.6]{q_0\mid\xix}$}; 
\draw [-] (f) -- (e);
\draw [->] (l) to [out=190,in=30] (f);

\end{tikzpicture}
\captionsetup{justification=centering,margin=2cm}
\captionsetup{justification=centering,margin=2cm}
\caption{Illustration of 1221.|411.}
\label{psi1_\ind21_pic3}
\end{center}
\end{figure}

\begin{minipage}[t]{0.5\textwidth}

 \begin{enumerate}[start=1222,leftmargin=20pt,label*={\arabic*}]
                   \item[1222.] If $q_0\nprec s_\tau$, then 
\end{enumerate}
\begin{multline*} \cph\phi_{\ind|2,1}(\vec{w}\ |\vec{q}\ ; z)=\\\cph=(\w_1,...,\w_\tau,q_0,s_{\tau+1},...,\w_\ind,z\ ;q_1)\in M^U_{\ind+2,1}.\end{multline*}

\end{minipage}
\vrule \hspace{0.5cm} 
\begin{minipage}[t]{0.5\textwidth}

\begin{enumerate}[start=1,leftmargin=15pt,label*={\arabic*}]

  \item[412.]  If ($\eta=\theta$ and $\xix\prec \wii_{\theta+1}$) or ($\eta=\theta+1$ and ($\wii_{\theta}\sim\wii_{\breve{\theta}+2}$ or $\wii_{\breve{\theta}}\nsucc\xix$ ) or $\eta>\breve{\theta}+1$, then
 $$\cps\psi^1_{\ind|2,1}(\vec{\wii}\ ;\xix) = (\wii_1,...,\wii_{\eta-1},\wii_{\eta+1},...,\wii_{\ind+1}\ |\wii_{\eta}, \xix\ ;\wii_{\ind+2}).$$

\end{enumerate}

\end{minipage}

\bigskip

\begin{figure}[H]
\begin{center}
\begin{tikzpicture} [>=latex,every node/.style={minimum width=3em, node distance=4em},scale=1.5,transform shape]

\node [draw,circle, fill=blue!20] (b) {$\Scale[0.6]{\w_1\mid\wii_1}$}; 
\node [node distance=2.5em,right of=b] (d) {$...$}; 
\node [draw,circle, fill=blue!20, node distance=2.5em,right of=d] (e) {$\Scale[0.6]{\w_{\tau}\mid\wii_{\eta-1}}$}; 

\node [node distance=3em,rotate=-35, above of=e, rotate=-90] (ee) {$\nprec$};
\node [node distance=1em, above right of=ee] (eee) {$\Scale[0.6]{q_1\mid\xix}$};
\node [node distance=2.5em,right of=e] (ef) {$\preceq$}; 
\node [draw,circle, fill=orange!20, node distance=5em,right of=e] (f) {$\Scale[0.6]{q_0\mid\wii_{\eta}}$}; 
 
\node [node distance=2.5em,right of=f] (ff) {$\preceq$}; 
\node [draw,circle, fill=blue!20, node distance=5em,right of=f] (fff) {$\Scale[0.6]{\w_{\tau+1}\mid\wii_{\eta+1}}$}; 
\node [node distance=2.5em,rotate=-55, above of=fff, rotate=-90] (ffff) {$\prec$}; 
\node [node distance=2.5em,right of=fff] (g) {$\prec$}; 
\node [node distance=1em,right of=g] (h) {$...$}; 
\node [node distance=1em,right of=h] (i) {$\prec$}; 
\node [draw,circle, fill=blue!20, node distance=2.5em,right of=i] (ii) {$\Scale[0.6]{\w_\ind|\wii_{\ind+1}}$}; 
\node [node distance=2.5em,right of=ii] (iii) {$\prec$}; 
\node [draw,circle, fill=orange!50, node distance=2.5em,right of=iii] (j)  {$\Scale[0.6]{z|\wii_{\ind+2}}$}; 
\node [node distance=3.25em, above of=j,rotate=-90] (k) {$\prec$}; 
\node [draw,circle, fill=orange!20, node distance=3.25em, above of=k] (l) {$\Scale[0.6]{q_1|\xix}$}; 

\draw [<-] (l) to [out=190,in=30] (f);

\end{tikzpicture}
\captionsetup{justification=centering,margin=2cm}
\caption{Illustration of 1222.|412.}
\label{psi1_\ind21_pic4}
\end{center}
\end{figure}

Here we have a special case when $\tau=\theta$. Note that it is possible that  $w_\tau\succ q_1$:

\begin{figure}[H]
\begin{center}
\begin{tikzpicture} [>=latex,every node/.style={minimum width=3em, node distance=4em},scale=1.5,transform shape]

\node [draw,circle, fill=blue!20] (b) {$\Scale[0.6]{\w_1\mid\wii_1}$};
\node [node distance=2.5em,right of=b] (d) {$...$}; 
\node [draw,circle, fill=blue!20, node distance=2.5em,right of=d] (e)  {$\Scale[0.6]{\w_{\theta}\mid\wii_{\eta-1}}$}; 
\node [draw,circle, fill=orange!20, node distance=5em,right of=e] (f) {$\Scale[0.6]{q_0\mid\wii_{\eta}}$}; 
\node [node distance=2.5em,right of=f] (ff) {$\preceq$}; 
\node [draw,circle, fill=blue!20, node distance=5em,right of=f] (fff) {$\Scale[0.6]{\w_{\theta+1}\mid\wii_{\eta+1}}$}; 
\node [node distance=2.5em,rotate=-55, above of=fff, rotate=-90] (ffff) {$\prec$}; 
\node [node distance=2.5em,right of=fff] (g) {$\prec$}; 
\node [node distance=1em,right of=g] (h) {$...$}; 
\node [node distance=1em,right of=h] (i) {$\prec$}; 
\node [draw,circle, fill=blue!20, node distance=2.5em,right of=i] (ii)  {$\Scale[0.6]{\w_\ind|\wii_{\ind+1}}$}; 
\node [node distance=2.5em,right of=ii] (iii) {$\prec$}; 
\node [draw,circle, fill=orange!50, node distance=2.5em,right of=iii] (j) {$\Scale[0.6]{z|\wii_{\ind+2}}$};  
\node [node distance=3.25em, above of=j,rotate=-90] (k) {$\prec$}; 
\node [draw,circle, fill=orange!20, node distance=3.25em, above of=k] (l)  {$\Scale[0.6]{q_1|\xix}$};  
\draw [-] (f) -- (e);
\draw [<-] (l) to [out=190,in=30] (f);
\draw [-] (e) to [out=-50,in=-130] (fff);

\end{tikzpicture}
\captionsetup{justification=centering,margin=2cm}
\caption{Illustration of 1222.|412.}
\label{psi1_\ind21_pic2}
\end{center}
\end{figure}

\begin{minipage}[t]{0.5\textwidth}

\begin{enumerate}[start=13,leftmargin=10pt,label*={\arabic*}]

\item [13.] If $z\prec \w_\ind$.

    \begin{enumerate}[leftmargin=5pt,label*={\arabic*}]
        \item [131.] If $\w_{\theta+1}\succ q_1$ and $\w_\theta\succ q_0$, then  $$\cph\phi_{\ind|2,1}(\vec{\w}\ |\vec{q}\ ;z)=(\vec{\w}\ ;\vec{q}\ ;z)\in M^U_{\ind,2,1}.$$

\end{enumerate}

\end{enumerate}

\end{minipage}
\vrule \hspace{0.5cm} 
\begin{minipage}[t]{0.5\textwidth}

\begin{enumerate}[start=1,leftmargin=5pt,label*={\arabic*}]

\item [5.] If $(\vec{\w}\ ;\kap_0,\kap_1;z) \in M^U_{\ind,2,1}$ and $z\succ\vec{\kap},$ and $z\prec\w_\ind,$ and $\w_{\theta+1}\succ q_1$ and $\w_\theta\succ q_0$, then

 $$\cps\psi^3_{\ind|2,1}(\vec{\w}\ ;\kap_0,\kap_1;z)=(\vec{\w}\ |\kap_0,\kap_1;z).$$

\end{enumerate}

\end{minipage}

\bigskip

\begin{figure}[H]
\begin{center}
\begin{tikzpicture} [>=latex,every node/.style={minimum width=3em, node distance=4em}]

\node [draw,circle, fill=blue!20] (b) {$\w_1$}; 
\node [node distance=2.5em,right of=b] (d) {$...$}; 
\node [draw,circle, fill=blue!20, node distance=2.5em,right of=d] (f) {$\w_{\theta}$}; 
\node [draw,circle, fill=blue!20, node distance=5em,right of=f] (ff) {$\w_{\theta+1}$}; 
\node [node distance=2.5em,right of=ff] (g) {$\prec$}; 
\node [node distance=1em,right of=g] (h) {$...$}; 
\node [node distance=1em,right of=h] (i) {$\prec$}; 
\node [draw,circle, fill=blue!20, node distance=2.5em,right of=i] (ii) {$\w_{\ind}$}; 
\node [node distance=2.5em, above of=ii,rotate=-90] (iii) {$\prec$}; 
\node [draw,circle, fill=orange!50, node distance=2.5em,above of=iii] (j) {$z$}; 
\node [node distance=2.5em,rotate=50, above of=j, rotate=-90] (m) {$\prec$}; 
\node [node distance=2.5em,rotate=70, above of=j, rotate=-90] (n) {$\prec$};

\node [node distance=5em, above of=f,rotate=-90] (k) {$\prec$}; 
\node [draw,circle, fill=orange!20, node distance=5em, above of=k] (l) {$q_0$}; 
\node [node distance=5em, above of=ff,rotate=-90] (kk) {$\prec$}; 
\node [draw,circle, fill=orange!20, node distance=5em, above of=kk] (ll) {$q_1$}; 
\draw [-] (f) -- (ff);
\draw [->] (l) -- (ll);

\end{tikzpicture}
\captionsetup{justification=centering,margin=2cm}
\caption{Illustration of 131.|5.}
\label{phi_\ind21_pic4a}
\end{center}
\end{figure}

\begin{minipage}[t]{0.5\textwidth}

\begin{enumerate}[start=132,leftmargin=15pt,label*={\arabic*}]

\item [132.] If $w_{\theta+1}\nsucc q_1$ or $w_\theta\nsucc q_0$, then  
   Let $$\gamma=\max(0,\{\theta<i<\ind| w_{i}\sim z\}).$$ 
    \begin{enumerate}[leftmargin=5pt,label*={\arabic*}]
       \item [1321.] if $\gamma>\theta$ and  $\gamma>\tau$ (i.e. $\exists\ w_\gamma\sim z$), then $$\cph\phi_{\ind|2,1}(\vec{w}\ |\vec{q}\ ;z)=(\vec{w}\ ;\vec{q}\ ;z)\in M^U_{\ind,2,1}.$$ 

This is the first type exceptional element, shown before on the Figure \ref{Firstn21} and on the picture below:

\end{enumerate}

\end{enumerate}

\end{minipage}
\vrule \hspace{0.5cm} 
\begin{minipage}[t]{0.5\textwidth}

\begin{enumerate}[start=1,leftmargin=5pt,label*={\arabic*}]

\item [6.] If $(\vec{w}\ ;\kap_0,\kap_1;z) \in M^U_{\ind,2,1}$  is the first type exceptional element , then
\\
\bigskip

 $$\cps\psi^3_{\ind|2,1}(\vec{w}\ ;\kap_0,\kap_1;z)=(\vec{w}\ |\kap_0,\kap_1;z).$$

\end{enumerate}

\end{minipage}

\bigskip

\begin{figure}[H]
\begin{center}
\begin{tikzpicture} [>=latex,every node/.style={minimum width=3em, node distance=4em}]

\node  (start) {$...$}; 
\node [draw,circle, fill=blue!20, node distance=2.5em, right of=start] (a)  {$\w_{\theta}$};
\node [draw,circle, fill=blue!20, right of=a] (b) {$\w_{\theta+1}$}; 
\node [node distance=2.5em,right of=b] (c) {$\prec$}; 
\node [node distance=1em,right of=c] (d) {$...$}; 
\node [node distance=1em,right of=d] (e) {$\prec$}; 
\node [draw,circle, fill=blue!20, node distance=2.5em,right of=e] (f) {$\w_{\gamma}$}; 
\node [node distance=2.5em,right of=f] (g) {$\prec$}; 
\node [node distance=1em,right of=g] (h) {$...$}; 
\node [node distance=1em,right of=h] (i) {$\prec$}; 
\node [draw,circle, fill=blue!20, node distance=2.5em,right of=i] (j) {$\w_{\ind}$}; 
\node [node distance=2.5em, above of=j,rotate=-90] (k) {$\prec$}; 
\node [draw,circle, fill=orange!50, node distance=2.5em,above of=k] (l) {$z$}; 
\node [node distance=3.5em, left of=l] (m) {$\prec$}; 
\node [draw,circle, fill=orange!20, node distance=3.5em,left of=m] (n) {$q_1$}; 
\node [node distance=2.5em, above of=f,rotate=-90] (o) {$\prec$}; 
\node [node distance=2.5em, above of=l,rotate=-90] (p) {$\prec$}; 
\node [draw,circle, fill=orange!20, node distance=2.5em,above of=p] (q) {$q_0$}; 
\draw [-] (b) -- (a);
\draw [-] (l) -- (f);
\draw [<-] (n) -- (q);

\end{tikzpicture}
\captionsetup{justification=centering,margin=2cm}
\caption{First exceptional element of $M^U_{\ind,2,1}$}
\label{Firstn21'}
\end{center}

\end{figure}

\begin{minipage}[t]{0.5\textwidth}

\begin{enumerate}[start=1322,leftmargin=20pt,label*={\arabic*}]

\item [1322.] Otherwise (i.e. if $\gamma=0$ or $\gamma<\tau$), we have: 
\\
\bigskip
\\
\bigskip
\\
\bigskip
\\
\bigskip
$$\cph\phi_{\ind|2,1}(\vec{w}\ |\vec{q}\ ;z)=(\phi_{\ind|2}(\vec{w}\ |\vec{q}\ );z).$$

\end{enumerate}

\end{minipage}
\vrule \hspace{0.5cm} 
\begin{minipage}[t]{0.5\textwidth}

\begin{enumerate}[start=1,leftmargin=5pt,label*={\arabic*}]

\item [7.]  For $(\vec{\wii}\ ;\xix) \in M^U_{\ind+2,1}$, such that $$\xix\prec\wii_{\ind+2} \text{ and } \eta=0,$$ what implies $$(\xix\succ\wii_{\theta+1} \text{ and }  \xix\succ\wii_{\theta}),$$  we have $$\cps\psi^2_{\ind|2,1}(\vec{\wii}\ ;\xix)=(\psi_{\ind|2}(\vec{\wii}\ );\xix).$$
\end{enumerate}

\end{minipage}

\bigskip

\begin{figure}[H]
\begin{center}
\begin{tikzpicture} [>=latex,every node/.style={minimum width=3em, node distance=4em},scale=1.5,transform shape]

\node [draw,circle, fill=blue!20] (b) {$\Scale[0.6]{\w_1\mid\wii_1}$};
\node [node distance=2.5em,right of=b] (d) {$...$}; 
\node [draw,circle, fill=orange!20, node distance=2.5em,right of=d] (f){$\Scale[0.6]{q_{01}\mid\wii_{\breve{\theta}}}$};
\node [node distance=2.5em,rotate=-54, above of=f, rotate=90] (mm1) {$\prec$}; 
\node [draw,circle, fill=orange!20, node distance=5em,right of=f] (fff) {$\Scale[0.6]{q_{01}\mid\wii_{\breve{\theta}+1}}$};
\node [node distance=2.5em,rotate=-56, above of=fff, rotate=90] (mm2) {$\prec$}; 
\node [node distance=2.5em,right of=fff] (g) {$\prec$}; 
\node [node distance=1em,right of=g] (h) {$...$}; 
\node [node distance=1em,right of=h] (i) {$\prec$}; 
\node [draw,circle, fill=blue!20, node distance=2.5em,right of=i] (ii) {$\Scale[0.6]{\w_{\ind-1}\mid\wii_{\ind+1}}$}; 
\node [node distance=2.5em,right of=ii] (iii) {$\prec$}; 
\node [draw,circle, fill=blue!20, node distance=2.5em,right of=iii] (j) {$\Scale[0.6]{\w_{\ind}\mid\wii_{\ind+2}}$};  
\node [node distance=2.5em, above of=j,rotate=-90] (k) {$\prec$}; 
\node [draw,circle, fill=orange!50, node distance=2.5em, above of=k] (l)  {$\Scale[0.6]{z\mid\xix}$}; 
\draw [<->] (f) -- (fff);

\end{tikzpicture}
\captionsetup{justification=centering,margin=2cm}
\caption{Illustration of 332.|19. (general case)}
\label{psi1_\ind21_pic6}
\end{center}
\end{figure}

\begin{minipage}[t]{0.5\textwidth}

\begin{enumerate}[start=2,leftmargin=5pt,label*={\arabic*}]

\item [2.] If $q_1\succ z$ and $q_0\sim z$.
  \begin{enumerate}[start=1,leftmargin=5pt,label*={\arabic*}]
    \item [21.] If $q_0\succ \vec{\w}$
    \begin{enumerate}[start=1,leftmargin=5pt,label*={\arabic*}]
       \item [211.] If $z\succ \vec{\w}$, then 

$$\cph\phi_{\ind|2,1}(\vec{\w}\ |\vec{q}\ ;z)=(\vec{\w}\ ;\vec{q}\ ;z)\in M^U_{\ind,2,1}.$$

\end{enumerate}
\end{enumerate}
\end{enumerate}
\end{minipage}
\vrule \hspace{0.5cm} 
\begin{minipage}[t]{0.5\textwidth}

\begin{enumerate}[start=1,leftmargin=5pt,label*={\arabic*}]

\item [8.] If   For $(\vec{\w}\ ;\kap_0,\kap_1;z) \in M^U_{\ind,2,1}$, such that
$$q_1\succ z\text{ and }q_0\sim z \text{ and } z\succ \vec{\w},$$
we have
$$\cps\psi^3_{\ind|2,1}(\vec{\w}\ ;\kap_0,\kap_1;z)=(\vec{\w}\ |\kap_0,\kap_1;z).$$

\end{enumerate}

\end{minipage}

\begin{figure}[H]
\begin{center}
\begin{tikzpicture} [>=latex,every node/.style={minimum width=3em, node distance=4em}]

\node [draw,circle, fill=blue!20] (b) {$\w_1$}; 
\node [node distance=2.5em,right of=b] (d) {$...$}; 
\node [draw,circle, fill=blue!20, node distance=2.5em,right of=d] (f) {$\w_\theta$}; 
\node [draw,circle, fill=blue!20, node distance=5em,right of=f] (ii) {$\w_{\theta+1}$}; 
\node [node distance=2.5em,right of=ii] (iii) {$\prec$};
\node [node distance=1em,right of=iii] (i1) {$...$};
 \node [node distance=1em,right of=i1] (i2) {$\prec$};
\node [draw,circle, fill=blue!20, node distance=2.5em,right of=i2] (j) {$\w_{\ind}$}; 
\node [node distance=2.5em, below of=j,rotate=-90] (k) {$\prec$}; 
\node [draw,circle, fill=orange!50, node distance=2.5em, below of=k] (l) {$z$}; 
\node [node distance=2.5em,rotate=25, above of=l, rotate=-90] (m) {$\prec$}; 
\node [node distance=2.5em,rotate=50, above of=l, rotate=-90] (m) {$\prec$}; 
\node [node distance=2.5em,rotate=75, above of=l, rotate=-90] (m) {$\prec$}; 

\node [node distance=2.5em, below of=l,rotate=-90] (q11) {$\prec$}; 
\node [draw,circle, fill=orange!20, node distance=2.5em, below of=q11] (q1) {$q_1$}; 
\node [draw,circle, fill=orange!20, node distance=5em, left of=q1] (q0) {$q_0$}; 
\draw [-] (f) -- (ii);
\draw [-] (l) -- (q0);
\draw [<-] (q1) -- (q0);

\end{tikzpicture}
\captionsetup{justification=centering,margin=2cm}
\caption{Illustration of 211.|8.}
\label{phi_\ind21_pic5}
\end{center}
\end{figure}
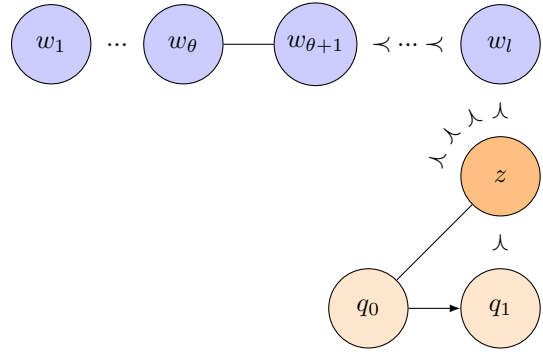

\begin{minipage}[t]{0.5\textwidth}

\begin{enumerate}[start=212,leftmargin=15pt,label*={\arabic*}]

       \item [212.] If $z\nsucc \vec{\w}$ and $z\nprec \w_\ind$, then  

$$\cph\phi_{\ind|2,1}(\vec{w}\ |\vec{q}\ ;z)=(\vec{w},z\ ;\vec{q}\ )\in M^U_{\ind+1,2}.$$

\end{enumerate}
\end{minipage}
\vrule \hspace{0.5cm} 
\begin{minipage}[t]{0.5\textwidth}

\begin{enumerate}[start=1,leftmargin=5pt,label*={\arabic*}]

\item [9.] For $(\vec{\w},\xix\ ; \kap_0,\kap_1) \in M^U_{\ind+1,2}$, such that 
$\kap_0\succ\vec{\w},$ $\kap_0\sim\xix$ and $\kap_1\succ \xix$, we have $$\cps\psi^2_{\ind|2,1}(\vec{\w},\xix\ ;\kap_0,\kap_1)=(\vec{\w}\ |\kap_0,\kap_1;\xix).$$

\end{enumerate}

\end{minipage}

\bigskip

\begin{figure}[H]
\begin{center}
\begin{tikzpicture} [>=latex,every node/.style={minimum width=3em, node distance=4em}]

\node [draw,circle, fill=blue!20] (b) {$\w_1$}; 
\node [node distance=2.5em,right of=b] (d) {$...$}; 
\node [draw,circle, fill=blue!20, node distance=2.5em,right of=d] (f) {$\w_\theta$}; 
\node [draw,circle, fill=blue!20, node distance=5em,right of=f] (ii) {$\w_{\theta+1}$}; 
\node [node distance=2.5em,right of=ii] (iii) {$\prec$};
\node [node distance=1em,right of=iii] (i1) {$...$};
 \node [node distance=1em,right of=i1] (i2) {$\prec$};
\node [draw,circle, fill=blue!20, node distance=2.5em,right of=i2] (j) {$\w_{\ind}$}; 
\node [node distance=2.5em, right of=j] (k) {$\preceq$}; 
\node [node distance=2.5em, below of=j,rotate=-90] (k) {$\prec$}; 
\node [draw,circle, fill=orange!50, node distance=5em, right of=j] (l) {$z\mid\xix$}; 

\node [node distance=2.5em, below of=l,rotate=-90] (q11) {$\prec$}; 
\node [draw,circle, fill=orange!20, node distance=2.5em, below of=q11] (q1) {$q_1$}; 
\node [draw,circle, fill=orange!20, node distance=5em, left of=q1] (q0) {$q_0$}; 
\node [node distance=2.5em,rotate=25, above of=q0, rotate=-90] (m) {$\prec$}; 
\node [node distance=2.5em,rotate=50, above of=q0, rotate=-90] (m) {$\prec$}; 
\node [node distance=2.5em,rotate=75, above of=q0, rotate=-90] (m) {$\prec$}; 
\draw [-] (f) -- (ii);
\draw [-] (l) -- (q0);
\draw [<-] (q1) -- (q0);

\end{tikzpicture}
\captionsetup{justification=centering,margin=2cm}
\caption{Illustration of 212.|9.}
\label{phi_\ind21_pic7}
\end{center}
\end{figure}
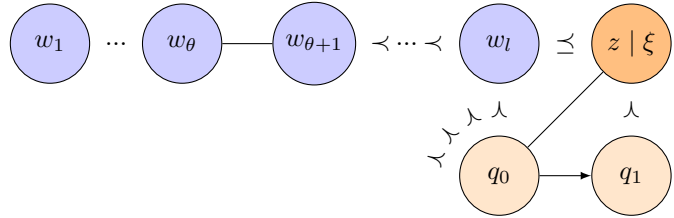

\begin{minipage}[t]{0.5\textwidth}

\begin{enumerate}[start=22,leftmargin=10pt,label*={\arabic*}]

    \item [22.] If $q_0\nsucc \vec{w}$ and $q_0\nprec \w_\ind$, then  

$$\cph\phi_{\ind|2,1}(\vec{w}\ |\vec{q}\ ;z)=(\vec{w},q_0,q_1\ ;z)\in M^U_{\ind+2,1}.$$

\end{enumerate}
\end{minipage}
\vrule \hspace{0.5cm} 
\begin{minipage}[t]{0.5\textwidth}

\begin{enumerate}[start=1,leftmargin=5pt,label*={\arabic*}]

\item [10.] For $(\vec{\wii}\ ;\xix) \in M^U_{\ind+2,1}$, such that \\$\wii_{\ind+1}\sim\xix$, $\wii_{\ind+1}\sim\wii_{\ind+2}$, and $\wii_{\ind+2}\succ\xix$, we have $$\cps\psi^2_{\ind|2,1}(\vec{\wii}\ ;\xix)=(\wii_1,...,\wii_\ind\I \wii_{\ind+1},\wii_{\ind+2}\;\xix).$$

\end{enumerate}

\end{minipage}
\bigskip

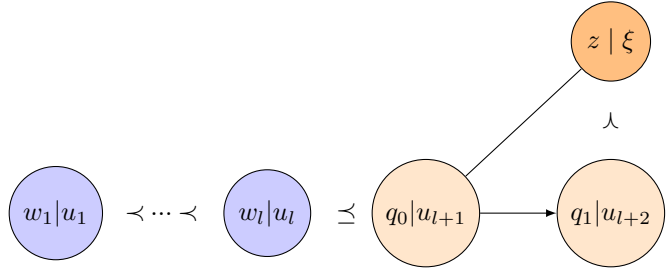
\begin{figure}[H]
\begin{center}
\begin{tikzpicture} [>=latex,every node/.style={minimum width=3em, node distance=4em}]

\node [draw,circle, fill=blue!20] (b) {$\w_1|\wii_1$}; 
 \node [node distance=3em,right of=b] (b1) {$\prec$};
\node [node distance=1em,right of=b1] (d) {$...$}; 
 \node [node distance=1em,right of=d] (i2) {$\prec$};
\node [draw,circle, fill=blue!20, node distance=3em,right of=i2] (j) {$\w_\ind|\wii_\ind$}; 
\node [node distance=3em, right of=j] (k) {$\preceq$}; 
\node [draw,circle, fill=orange!20, node distance=6em, right of=j] (q0) {$q_0|\wii_{\ind+1}$}; 
\node [draw,circle, fill=orange!20, node distance=7em, right of=q0] (q1) {$q_1|\wii_{\ind+2}$}; 
\node [node distance=3.5em, above of=q1,rotate=-90] (q11) {$\prec$}; 
\node [draw,circle, fill=orange!50, node distance=3em, above of=q11] (z) {$z\mid\xix$}; 
\draw [-] (z) -- (q0);
\draw [<-] (q1) -- (q0);
\end{tikzpicture}
\captionsetup{justification=centering,margin=2cm}
\caption{Illustration of 22.|10.}
\label{phi_\ind21_pic8}
\end{center}
\end{figure}

\begin{minipage}[t]{0.5\textwidth}

\begin{enumerate}[start=23,leftmargin=10pt,label*={\arabic*}]

            \item[23.] If $q_0\prec \w_\ind$.
            \begin{enumerate}[start=1,leftmargin=5pt,label*={\arabic*}]
                \item[231.] If $q_1\sim \w_\ind$, then
\\
\\
\\
 $$\cph\phi_{\ind|2,1}(\vec{\w}\ |\vec{q}\ ;z)=(\vec{\w},q_1;q_0,z)\in M^U_{\ind+1,2}.$$
\end{enumerate}

\end{enumerate}
\end{minipage}
\vrule \hspace{0.5cm} 
\begin{minipage}[t]{0.5\textwidth}

\begin{enumerate}[start=1,leftmargin=5pt,label*={\arabic*}]

\item [11.] For $(\vec{\w},\xix\ ;\kap_0,\kap_1) \in M^U_{\ind+1,2}$, such that $$\xix\sim \w_\ind \text{ and }\xix\sim \kap_0\text{ and }\kap_0\prec \w_\ind\text{ and }\kap_1\prec \xix,$$
we have $$\cps\psi^2_{\ind|2,1}(\vec{\w},\xix\ ;\kap_0,\kap_1)=(\vec{\w}\ |\kap_0,\xix\ ;\kap_1).$$

\end{enumerate}

\end{minipage}
\bigskip

\begin{figure}[H]
\begin{center}
\begin{tikzpicture} [>=latex,every node/.style={minimum width=3em, node distance=4em}]

\node [draw,circle, fill=blue!20] (b) {$\w_1$}; 
\node [node distance=2.5em,right of=b] (d) {$...$}; 
\node [draw,circle, fill=blue!20, node distance=2.5em,right of=d] (f) {$\w_\ind$}; 
\node [draw,circle, fill=orange!20, node distance=5em,right of=f] (j) {$q_1\mid\xix$};
\node [node distance=2.5em, above of=f,rotate=-90] (k) {$\prec$}; 
\node [draw,circle, fill=orange!20, node distance=2.5em, above of=k] (l) {$q_0$}; 
\node [node distance=2.5em, above of=j,rotate=-90] (o) {$\prec$};
\node [draw,circle, fill=orange!50, node distance=2.5em, above of=o] (p) {$z\mid q_1$};
\draw [-] (l) -- (p);
\draw [-] (f) -- (j);
\draw [->] (l) -- (j);

\end{tikzpicture}
\captionsetup{justification=centering,margin=2cm}
\caption{Illustration of 231.|11.}
\label{phi_\ind21_pic8}
\end{center}
\end{figure}
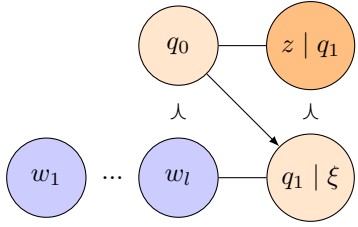

\begin{minipage}[t]{0.5\textwidth}

\begin{enumerate}[start=232,leftmargin=15pt,label*={\arabic*}]

       \item [232.] If $q_1\prec \w_\ind$ and $\w_{\theta}\succ q_0$ and $\w_{\theta+1}\succ q_1$, then 
\\
\\
\\
$$\cph\phi_{\ind|2,1}(\vec{\w}\ |\vec{q}\ ;z)=(\vec{\w}\ ;\vec{q}\ ;z)\in M^U_{\ind,2,1}. $$

\end{enumerate}
\end{minipage}
\vrule \hspace{0.5cm} 
\begin{minipage}[t]{0.5\textwidth}

\begin{enumerate}[start=1,leftmargin=5pt,label*={\arabic*}]

\item [12.] For $(\vec{\w}\ ;\kap_0,\kap_1;z) \in M^U_{\ind,2,1}$, such that 
\\$\w_\theta\succ\kap_0 \text{ and }\w_{\theta+1}\succ\kap_1\text{ and }\kap_0\sim z\text{ and }\kap_1\succ z $, we have
 $$\cps\psi^3_{\ind|2,1}(\vec{\w}\ ;\kap_0,\kap_1;z)=(\vec{\w}\ |\kap_0,\kap_1;z).$$

\end{enumerate}

\end{minipage}
\bigskip

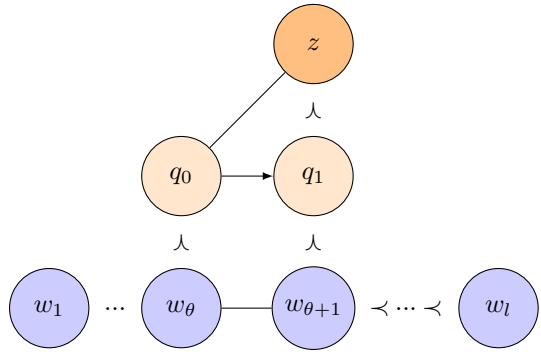
\begin{figure}[H]
\begin{center}
\begin{tikzpicture} [>=latex,every node/.style={minimum width=3em, node distance=4em}]

\node [draw,circle, fill=blue!20] (b) {$\w_1$}; 
\node [node distance=2.5em,right of=b] (d) {$...$}; 
\node [draw,circle, fill=blue!20, node distance=2.5em,right of=d] (f) {$\w_{\theta}$}; 
\node [draw,circle, fill=blue!20, node distance=5em,right of=f] (ff) {$\w_{\theta+1}$}; 
\node [node distance=2.5em,right of=ff] (g) {$\prec$}; 
\node [node distance=1em,right of=g] (h) {$...$}; 
\node [node distance=1em,right of=h] (i) {$\prec$}; 
\node [draw,circle, fill=blue!20, node distance=2.5em,right of=i] (ii) {$\w_{\ind}$};

\node [node distance=2.5em, above of=f,rotate=-90] (k) {$\prec$}; 
\node [draw,circle, fill=orange!20, node distance=2.5em, above of=k] (l) {$q_0$}; 
\node [node distance=2.5em, above of=ff,rotate=-90] (kk) {$\prec$}; 
\node [draw,circle, fill=orange!20, node distance=2.5em, above of=kk] (ll) {$q_1$}; 

\node [node distance=2.5em, above of=ll,rotate=-90] (lla) {$\prec$}; 
\node [draw,circle, fill=orange!50, node distance=2.5em,above of=lla] (j) {$z$};

\draw [-] (f) -- (ff);
\draw [-] (l) -- (j);
\draw [->] (l) -- (ll);

\end{tikzpicture}
\captionsetup{justification=centering,margin=2cm}
\caption{Illustration of 232.|12.}
\label{phi_\ind21_pic321}
\end{center}
\end{figure}

\begin{minipage}[t]{0.5\textwidth}

\begin{enumerate}[start=233,leftmargin=15pt,label*={\arabic*}]

        \item [233.] If $q_1\prec \w_\ind$ and ($\w_{\theta}\nsucc q_0$ or $\w_{\theta+1}\nsucc q_1$), let $$\tau=\max(\{i<\ind| \w_{i}\nsucc q_1\vee \w_i\sim \w_{i+1}\}).$$

\begin{enumerate}[start=1,leftmargin=5pt,label*={\arabic*}]

                   \item[2331.] If $q_0\prec \w_\tau$ ($\Rightarrow\ q_1\sim \w_\tau$), then $$\cph\phi_{\ind|2,1}(\vec{w}\ |\vec{q}\ ;z)=(\vec{\w}\ ;\vec{q}\ ;z)\in M^U_{\ind,2,1}.$$

This case is isomorphic to the second exceptional element type, shown below:
\end{enumerate}

\end{enumerate}
\end{minipage}
\vrule \hspace{0.5cm} 
\begin{minipage}[t]{0.5\textwidth}

\begin{enumerate}[start=1,leftmargin=5pt,label*={\arabic*}]

\item [13.] If $(\vec{\w}\ ;\kap_0,\kap_1;z) \in M^U_{\ind,2,1}$ has the second\\ exceptional element type, then
\\
\\

 $$\cps\psi^3_{\ind|2,1}(\vec{\w}\ ;\kap_0,\kap_1;z)=(\vec{\w}\ |\kap_0,\kap_1;z).$$ 

\end{enumerate}

\end{minipage}
\bigskip

\begin{figure}[H]
\begin{center}
\begin{tikzpicture} [>=latex,every node/.style={minimum width=3em, node distance=4em}]

\node  (start) {$...$}; 
\node [draw,circle, fill=blue!20, node distance=2.5em, right of=start] (a)  {$\w_{\theta}$};
\node [draw,circle, fill=blue!20, right of=a] (b) {$\w_{\theta+1}$}; 
\node [node distance=2.5em,right of=b] (c) {$\prec$}; 
\node [node distance=1em,right of=c] (d) {$...$}; 
\node [node distance=1em,right of=d] (e) {$\prec$}; 
\node [draw,circle, fill=blue!20, node distance=2.5em,right of=e] (f) {$\w_{\tau}$}; 
\node [node distance=2.5em,right of=f] (g) {$\prec$}; 
\node [node distance=1em,right of=g] (h) {$...$}; 
\node [node distance=1em,right of=h] (i) {$\prec$}; 
\node [draw,circle, fill=blue!20, node distance=2.5em,right of=i] (j) {$\w_\ind$}; 
\node [node distance=2.5em, above of=j,rotate=-90] (k) {$\prec$}; 
\node [draw,circle, fill=orange!20, node distance=2.5em, above of=k] (l) {$q_1$}; 

\node [draw,circle, fill=orange!20, node distance=7em,left of=l] (n) {$q_0$}; 
\node [node distance=2.5em, above of=f,rotate=-90] (o) {$\prec$}; 
\node [node distance=2.5em, above of=l,rotate=-90] (p) {$\prec$}; 
\node [draw,circle, fill=orange!50, node distance=2.5em,above of=p] (q) {$z$}; 
\draw [-] (b) -- (a);
\draw [-] (l) -- (f);
\draw [<-] (l) -- (n);
\draw [-] (n) -- (q);

\end{tikzpicture}
\captionsetup{justification=centering,margin=2cm}
\caption{Illustration of 2331.|13. \\(Second exceptional element from $M^U_{\ind,2,1}$\ ).}
\label{Secondn21a}
\end{center}
\end{figure}
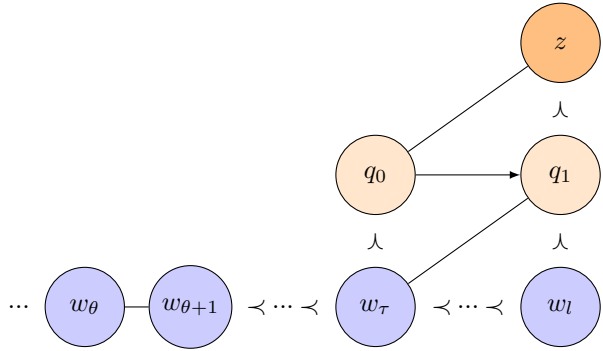

\begin{minipage}[t]{0.5\textwidth}

\begin{enumerate}[start=2332,leftmargin=20pt,label*={\arabic*}]
                    \item [2332.] If $q_0\nprec w_\tau$, then
\end{enumerate}
 \begin{multline*}\cph\phi_{\ind|2,1}(\vec{w}\ |\vec{q}\ ;z)=(\phi_{\ind|2}(\vec{w}\ |q_0,z)\ ;q_1)=\\\cph=(\w_1,...,w_\tau,q_0,z,w_{\tau+1},...,\w_\ind\ ;q_1)\in M^U_{\ind+2,1}.\end{multline*} The following 3 pictures illustrate this case:

\end{minipage}
\vrule \hspace{0.5cm} 
\begin{minipage}[t]{0.5\textwidth}

\begin{enumerate}[start=1,leftmargin=5pt,label*={\arabic*}]
                   \item [14.] For $(\vec{\wii}\ ;\xix) \in M^U_{\ind+2,1}$, such that $\xix\prec\wii_{\ind+2}$,  $\eta=\theta$ and $\xix\succ\wii_{\theta+1}$, we have $$\cps\psi^1_{\ind|2,1}(\vec{\wii}\ ;\xix)=(\wii_1,...,\wii_{\theta-1},\wii_{\theta+2},...,\wii_{\ind+2}\ |\wii_{\theta},\xix\ ;\wii_{\theta+1}).$$
As a reminder, \begin{align*} & \breve{\theta}=\max(\{i<\ind+2|\ \wii_i\sim \wii_{i+1}\}), \\ &   \eta=\max(0,\{i\geq\theta| \wii_i\sim \xix\}).\end{align*}

\end{enumerate}

\end{minipage}
\bigskip

\begin{figure}[H]
\begin{center}
\begin{tikzpicture} [>=latex,every node/.style={minimum width=3em, node distance=4em},scale=1.5,transform shape]

\node [draw,circle, fill=blue!20] (b) {$\Scale[0.6]{\w_1\mid\wii_1}$};
\node [node distance=2.5em,right of=b] (d) {$...$}; 

\node [draw,circle, fill=blue!20, node distance=2.5em,right of=d] (dd) {$\Scale[0.6]{\w_{\theta}\mid\wii_{\breve{\theta}-1}}$};
\node [draw,circle, fill=orange!20, node distance=5em,right of=dd] (f) {$\Scale[0.6]{q_0\mid\wii_{\breve{\theta}}}$};
\node [draw,circle, fill=orange!50, node distance=5em,right of=f] (fff)  {$\Scale[0.6]{z\mid\wii_{\breve{\theta}+1}}$};
\node [node distance=2.5em,rotate=-55, above of=fff, rotate=-90] (ffff) {$\succ$}; 
\node [node distance=2.5em,right of=fff] (g) {$\prec$}; 
\node [draw,circle, fill=blue!20, node distance=2.5em,right of=g] (t1)  {$\Scale[0.6]{\w_{\theta+1}\mid\wii_{\breve{\theta}+2}}$};
\node [node distance=2.5em,right of=t1] (t2) {$\prec$}; 
\node [node distance=1em,right of=t2] (h) {$...$}; 
\node [node distance=1em,right of=h] (i) {$\prec$}; 
\node [draw,circle, fill=blue!20, node distance=2.5em,right of=i] (j) {$\Scale[0.6]{\w_\ind\mid\wii_{\ind+2}}$}; 
\node [node distance=3.25em, above of=j,rotate=-90] (k) {$\prec$}; 
\node [draw,circle, fill=orange!20, node distance=3.25em, above of=k] (l)  {$\Scale[0.6]{q_1\mid\xix}$};
\draw [-] (f) -- (fff);
\draw [<-] (l) -- (f);
\draw [-] (dd) -- (f);
\draw [-] (dd) to [out=-50,in=-130] (t1);
\end{tikzpicture}
\captionsetup{justification=centering,margin=2cm}
\caption{Illustration of 2332.|14.\\ ($\phi_{\ind|2,1}$, if $\tau=\theta$. Note that it may happen that $w_\tau\succ q_1$).}
\label{phi_\ind21_pic9}
\end{center}
\end{figure}
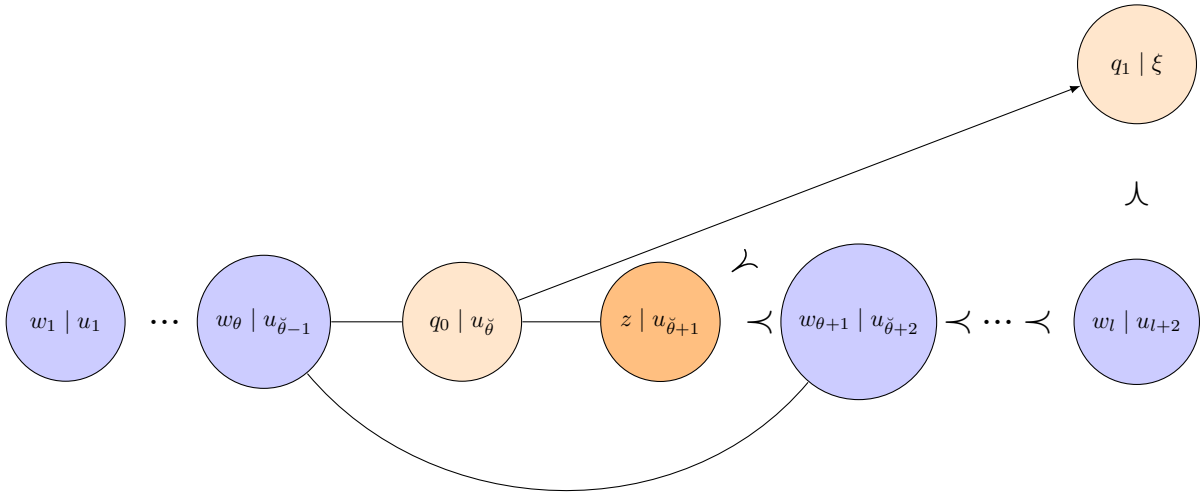

\begin{figure}[H]
\begin{center}
\begin{tikzpicture} [>=latex,every node/.style={minimum width=3em, node distance=4em},scale=1.5,transform shape]

\node [draw,circle, fill=blue!20] (b)  {$\Scale[0.6]{\w_1\mid\wii_1}$};
\node [node distance=2.5em,right of=b] (d) {$...$}; 

\node [draw,circle, fill=blue!20, node distance=2.5em,right of=d] (dd) {$\Scale[0.6]{\w_{\tau}\mid\wii_{\breve{\theta}-1}}$};
\node [draw,circle, fill=orange!20, node distance=5em,right of=dd] (f)  {$\Scale[0.6]{q_0\mid\wii_{\breve{\theta}}}$};
\node [draw,circle, fill=orange!50, node distance=5em,right of=f] (fff) {$\Scale[0.6]{z\mid\wii_{\breve{\theta}+1}}$};
\node [node distance=2.5em,rotate=-55, above of=fff, rotate=-90] (ffff) {$\succ$}; 
\node [node distance=2.5em,right of=fff] (g) {$\prec$}; 
\node [draw,circle, fill=blue!20, node distance=2.5em,right of=g] (t1) {$\Scale[0.6]{\w_{\tau+1}\mid\wii_{\breve{\theta}+2}}$};
\node [node distance=2.5em,rotate=-34, above of=t1, rotate=-90] (tb) {$\prec$}; 
\node [node distance=2.5em,right of=t1] (t2) {$\prec$}; 
\node [node distance=1em,right of=t2] (h) {$...$}; 
\node [node distance=1em,right of=h] (i) {$\prec$}; 
\node [draw,circle, fill=blue!20, node distance=2.5em,right of=i] (j) {$\Scale[0.6]{\w_\ind\mid\wii_{\ind+2}}$}; 
\node [node distance=4em, above of=j,rotate=-90] (k) {$\prec$}; 
\node [draw,circle, fill=orange!20, node distance=4em, above of=k] (l)  {$\Scale[0.6]{q_1\mid\xix}$};
\draw [-] (f) -- (fff);
\draw [<-] (l) -- (f);
\draw [-] (dd) -- (f);
\draw [-] (dd) -- (l);
\end{tikzpicture}
\captionsetup{justification=centering,margin=2cm}
\caption{Illustration of 2332.|14.\\ ($\phi_{\ind|2,1}$, if $q_1\sim w_\tau$ and $w_\tau\sim q_0$.)}
\label{phi_\ind21_pic10}
\end{center}
\end{figure}
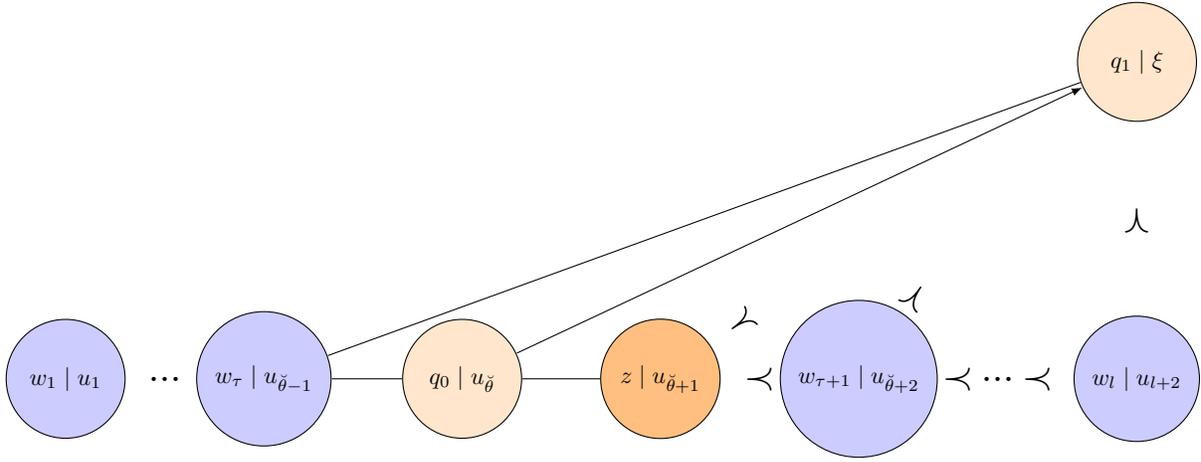

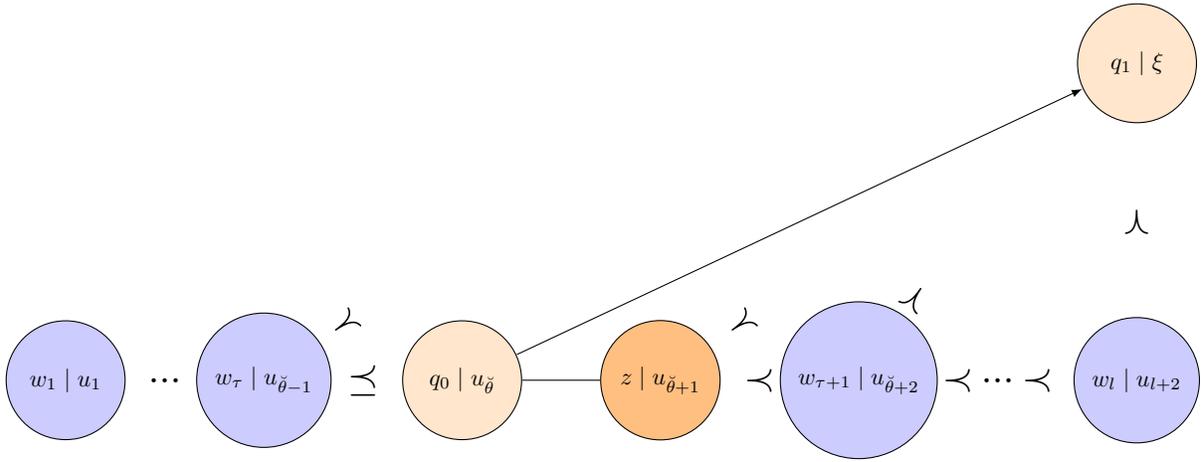
\begin{figure}[H]
\begin{center}
\begin{tikzpicture} [>=latex,every node/.style={minimum width=3em, node distance=4em},scale=1.5,transform shape]

\node [draw,circle, fill=blue!20] (b)  {$\Scale[0.6]{\w_1\mid\wii_1}$};
\node [node distance=2.5em,right of=b] (d) {$...$}; 

\node [draw,circle, fill=blue!20, node distance=2.5em,right of=d] (dd) {$\Scale[0.6]{\w_{\tau}\mid\wii_{\breve{\theta}-1}}$};
\node [node distance=2.5em,rotate=-55, above of=dd, rotate=-90] (dd1) {$\succ$}; 
 \node [node distance=2.5em,right of=dd] (d1) {$\preceq$}; 
\node [draw,circle, fill=orange!20, node distance=5em,right of=dd] (f) {$\Scale[0.6]{q_0\mid\wii_{\breve{\theta}}}$};
\node [draw,circle, fill=orange!50, node distance=5em,right of=f] (fff) {$\Scale[0.6]{z\mid\wii_{\breve{\theta}+1}}$};
\node [node distance=2.5em,rotate=-55, above of=fff, rotate=-90] (ffff) {$\succ$}; 
\node [node distance=2.5em,right of=fff] (g) {$\prec$}; 
\node [draw,circle, fill=blue!20, node distance=2.5em,right of=g] (t1){$\Scale[0.6]{\w_{\tau+1}\mid\wii_{\breve{\theta}+2}}$};
\node [node distance=2.5em,rotate=-34, above of=t1, rotate=-90] (tb) {$\prec$}; 
\node [node distance=2.5em,right of=t1] (t2) {$\prec$}; 
\node [node distance=1em,right of=t2] (h) {$...$}; 
\node [node distance=1em,right of=h] (i) {$\prec$}; 
\node [draw,circle, fill=blue!20, node distance=2.5em,right of=i] (j) {$\Scale[0.6]{\w_\ind\mid\wii_{\ind+2}}$};  
\node [node distance=4em, above of=j,rotate=-90] (k) {$\prec$}; 
\node [draw,circle, fill=orange!20, node distance=4em, above of=k] (l) {$\Scale[0.6]{q_1\mid\xix}$};
\draw [-] (f) -- (fff);
\draw [<-] (l) -- (f);
\end{tikzpicture}
\captionsetup{justification=centering,margin=2cm}
\caption{Illustration of 2332.|14.\\ ($\phi_{\ind|2,1}$, if $q_1\succ w_\tau$.)}
\label{phi_\ind21_pic11}
\end{center}
\end{figure}

\begin{minipage}[t]{0.5\textwidth}

\begin{enumerate}[start=3,leftmargin=5pt,label*={\arabic*}]

       \item [3.] If $q_1\succ z$ and $q_0\succ z$. 

This is the easiest case, since we insert $\vec{q}$ and $z$ independently. 
\begin{enumerate}[start=1,leftmargin=5pt,label*={\arabic*}]
   \item [31.] If $z\succ\vec{\w}$, then $$\cph\phi_{\ind|2,1}(\vec{\w}\ |\vec{q}\ ;z)=(\vec{\w}\ |\vec{q}\ ;z)\in M^U_{\ind,2,1}.$$

\end{enumerate} 
\end{enumerate}
\end{minipage}
\vrule \hspace{0.5cm} 
\begin{minipage}[t]{0.5\textwidth}

\begin{enumerate}[start=1,leftmargin=5pt,label*={\arabic*}]

\item [15.] For $(\vec{\w}\ ;\kap_0,\kap_1;z) \in M^U_{\ind,2,1}$, such that

 $z\succ\vec{\w}$, $\kap_0\succ z$ and $\kap_1\succ z$, 

we have
\bigskip
 $$\cps\psi^3_{\ind|2,1}(\vec{\w}\ ;\kap_0,\kap_1;z)=(\vec{\w}\ |\kap_0,\kap_1;z).$$ 

\end{enumerate}

\end{minipage}

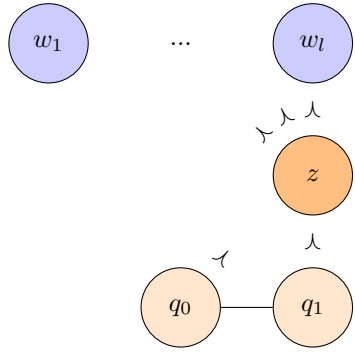
\begin{figure}[H]
\begin{center}
\begin{tikzpicture} [>=latex,every node/.style={minimum width=3em, node distance=4em}]

\node [draw,circle, fill=blue!20] (b) {$\w_1$}; 
\node [node distance=5em,right of=b] (d) {$...$}; 
\node [draw,circle, fill=blue!20, node distance=5em,right of=d] (f) {$\w_{\ind}$}; 
\node [node distance=2.5em, below of=f,rotate=-90] (k) {$\prec$}; 
\node [draw,circle, fill=orange!50, node distance=2.5em, below of=k] (l) {$z$};
\node [node distance=2.5em, below of=l,rotate=-90] (zl) {$\prec$}; 
\node [draw,circle, fill=orange!20, node distance=2.5em, below of=zl] (zq1) {$q_1$};
\node [draw,circle, fill=orange!20, node distance=5em, left of=zq1] (zq0) {$q_0$};

\node [node distance=2.5em,rotate=-40, above of=zq0, rotate=90] (mm) {$\succ$};  
\node [node distance=2.5em,rotate=25, above of=l, rotate=-90] (m) {$\prec$}; 
\node [node distance=2.5em,rotate=50, above of=l, rotate=-90] (n) {$\prec$};

\draw [-] (zq1) -- (zq0);
\end{tikzpicture}
\captionsetup{justification=centering,margin=2cm}
\caption{Illustration of 31.|15.}
\label{phi_\ind21_pic11a}
\end{center}
\end{figure}

\begin{minipage}[t]{0.5\textwidth}

\begin{enumerate}[start=32,leftmargin=10pt,label*={\arabic*}]

   \item [32.] If $z\nprec \w_\ind$ and $z\nsucc\vec{w}$ ($\Rightarrow q_0\nprec \w_\ind$)
      \begin{enumerate}[start=1,leftmargin=5pt,label*={\arabic*}]
          \item [321.] If $q_0\succ\vec{w}$, then $$\cph\phi_{\ind|2,1}(\vec{w}\ |\vec{q}\ ;z)=(\vec{w},z\ ;\vec{q}\ )\in M^U_{\ind+1,2}.$$
\end{enumerate}

\end{enumerate}
\end{minipage}
\vrule \hspace{0.5cm} 
\begin{minipage}[t]{0.5\textwidth}

\begin{enumerate}[start=1,leftmargin=5pt,label*={\arabic*}]

\item [16.]  If $(\vec{\w},\xix\ ;\kap_0,\kap_1) \in M^U_{\ind+1,2}$, \\and $\kap_0\succ\vec{\w}$, and $\kap_1\succ \xix$, then 

$$\cps\psi^2_{\ind|2,1}(\vec{\w},\xix\ ;\kap_0,\kap_1)=(\vec{\w}\ |\kap_0,\kap_1;\xix).$$

\end{enumerate}

\end{minipage}
\bigskip

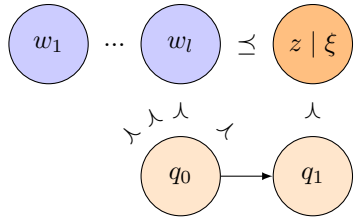
\begin{figure}[H]
\begin{center}
\begin{tikzpicture} [>=latex,every node/.style={minimum width=3em, node distance=4em}]

\node [draw,circle, fill=blue!20] (b) {$\w_1$}; 
\node [node distance=2.5em,right of=b] (d) {$...$}; 
\node [draw,circle, fill=blue!20, node distance=2.5em,right of=d] (f) {$\w_\ind$}; 
\node [node distance=2.5em,right of=f] (i) {$\preceq$}; 
\node [draw,circle, fill=orange!50, node distance=2.5em,right of=i] (j) {$z\mid\xix$};
\node [node distance=2.5em, below of=f,rotate=-90] (k) {$\prec$}; 
\node [draw,circle, fill=orange!20, node distance=2.5em, below of=k] (l) {$q_0$};
\node [node distance=2.5em,rotate=-45, above of=l, rotate=90] (mm) {$\succ$};  
\node [node distance=2.5em,rotate=25, above of=l, rotate=-90] (m) {$\prec$}; 
\node [node distance=2.5em,rotate=50, above of=l, rotate=-90] (n) {$\prec$}; 
\node [node distance=2.5em, below of=j,rotate=-90] (o) {$\prec$}; 
\node [draw,circle, fill=orange!20, node distance=2.5em, below of=o] (p) {$q_1$};
\draw [->] (l) -- (p);

\end{tikzpicture}
\captionsetup{justification=centering,margin=2cm}
\caption{Illustration of 321.|16.}
\label{phi_\ind21_pic12}
\end{center}
\end{figure}

\begin{minipage}[t]{0.5\textwidth}

\begin{enumerate}[start=322,leftmargin=15pt,label*={\arabic*}]

          \item [322.] If $q_0\nsucc\vec{\w}$,  then
\\
\\
\\
\\
\\
\\

 $$\cph\phi_{\ind|2,1}(\vec{\w}\ |\vec{q}\ ;z)=(\vec{\w}\ ,\vec{q}\ ;z)\in M^U_{\ind+2,1}.$$

\end{enumerate}

\end{minipage}
\vrule \hspace{0.5cm} 
\begin{minipage}[t]{0.5\textwidth}

\begin{enumerate}[start=1,leftmargin=5pt,label*={\arabic*}]

\item  [17.] For $(\vec{\wii}\ ;\xix) \in M^U_{\ind+2,1}$, such that $$\xix\prec\wii_{\ind+2}\text{ and }\xix\prec\wii_{\ind+2},$$ what implies $$ (\xix\prec\wii_{\theta+1} \text{ and  } \xix\prec\wii_{\theta}),$$  we have $$\cps\psi^2_{\ind|2,1}(\vec{\wii}\ ;\xix)=(\psi_{\ind|2}(\vec{\wii}\ );\xix)=(\wii_1,...,\wii_\ind| \wii_{\ind+1},\wii_{\ind+2};\xix).$$
\end{enumerate}

\end{minipage}

\begin{figure}[H]
\begin{center}
\begin{tikzpicture} [>=latex,every node/.style={minimum width=3em, node distance=4em}]

\node [draw,circle, fill=blue!20] (b) {$\w_1|\wii_1$}; 
 \node [node distance=3em,right of=b] (b1) {$\prec$};
\node [node distance=1em,right of=b1] (d) {$...$}; 
 \node [node distance=1em,right of=d] (i2) {$\prec$};
\node [draw,circle, fill=blue!20, node distance=3em,right of=i2] (j) {$\w_\ind|\wii_\ind$}; 
\node [node distance=3em, right of=j] (k) {$\preceq$}; 
\node [draw,circle, fill=orange!20, node distance=6em, right of=j] (q0) {$q_0|\wii_{\ind+1}$}; 
\node [node distance=3.5em,rotate=-40, above of=q0, rotate=90] (mm) {$\succ$}; 
\node [draw,circle, fill=orange!20, node distance=7em, right of=q0] (q1) {$q_1|\wii_{\ind+2}$}; 
\node [node distance=3.5em, above of=q1,rotate=-90] (q11) {$\prec$}; 
\node [draw,circle, fill=orange!50, node distance=3em, above of=q11] (z) {$z\mid\xix$}; 

\draw [<-] (q1) -- (q0);
\end{tikzpicture}
\captionsetup{justification=centering,margin=2cm}
\caption{Illustration of 322.|17.}
\label{phi_\ind21_pic8}
\end{center}
\end{figure}
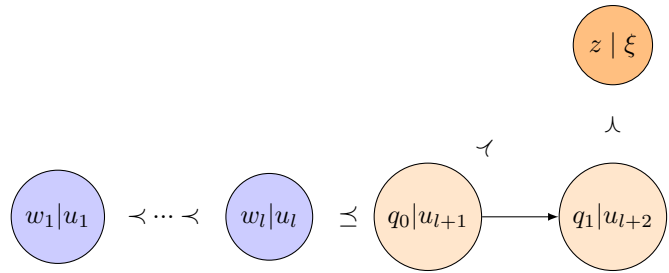

\begin{minipage}[t]{0.5\textwidth}

\begin{enumerate}[start=33,leftmargin=10pt,label*={\arabic*}]
    
   \item [33.] If $z\prec \w_\ind$
\begin{enumerate}[start=1,leftmargin=5pt,label*={\arabic*}]
             \item [331.] If $(\vec{\w}\ ;\vec{q}\ )\in M_{\ind,2}^U$, then 
$$\cph\phi_{\ind|2,1}(\vec{\w}\ |\vec{q}\ ;z)=(\vec{\w}\ ;\vec{q}\ ;z)\in M^U_{\ind,2,1}.$$

         \end{enumerate}

\end{enumerate}

\end{minipage}
\vrule \hspace{0.5cm} 
\begin{minipage}[t]{0.5\textwidth}

\begin{enumerate}[start=1,leftmargin=5pt,label*={\arabic*}]

\item [18.] If $(\vec{\w}\ ;\kap_0,\kap_1;z) \in M^U_{\ind,2,1}$, and 
$$\kap_0\succ z,\ \kap_1\succ z\text{ and } w_\ind\succ z, \text{ then}$$
$$\cps\psi^3_{\ind|2,1}(\vec{\w}\ ;\vec{\kap}\ ;z)=(\vec{\w}\ |\vec{\kap}\ ;z).$$ 

\end{enumerate}

\end{minipage}
\bigskip

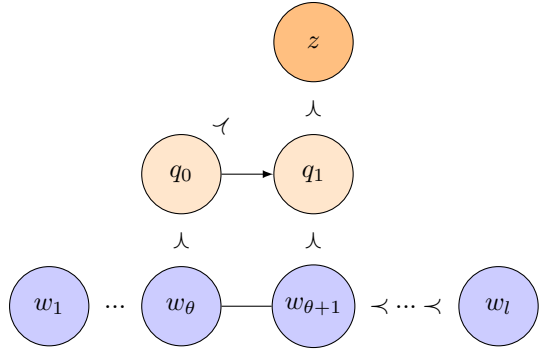
\begin{figure}[H]
\begin{center}
\begin{tikzpicture} [>=latex,every node/.style={minimum width=3em, node distance=4em}]

\node [draw,circle, fill=blue!20] (b) {$\w_1$}; 
\node [node distance=2.5em,right of=b] (d) {$...$}; 
\node [draw,circle, fill=blue!20, node distance=2.5em,right of=d] (f) {$\w_{\theta}$}; 
\node [draw,circle, fill=blue!20, node distance=5em,right of=f] (ff) {$\w_{\theta+1}$}; 
\node [node distance=2.5em,right of=ff] (g) {$\prec$}; 
\node [node distance=1em,right of=g] (h) {$...$}; 
\node [node distance=1em,right of=h] (i) {$\prec$}; 
\node [draw,circle, fill=blue!20, node distance=2.5em,right of=i] (ii) {$\w_{\ind}$};

\node [node distance=2.5em, above of=f,rotate=-90] (k) {$\prec$}; 
\node [draw,circle, fill=orange!20, node distance=2.5em, above of=k] (l) {$q_0$}; 
\node [node distance=2.5em, above of=ff,rotate=-90] (kk) {$\prec$}; 
\node [node distance=2.5em,rotate=-40, above of=l, rotate=90] (mm) {$\succ$};  
\node [draw,circle, fill=orange!20, node distance=2.5em, above of=kk] (ll) {$q_1$}; 

\node [node distance=2.5em, above of=ll,rotate=-90] (lla) {$\prec$}; 
\node [draw,circle, fill=orange!50, node distance=2.5em,above of=lla] (j) {$z$};

\draw [-] (f) -- (ff);
\draw [->] (l) -- (ll);

\end{tikzpicture}
\captionsetup{justification=centering,margin=2cm}
\caption{Illustration of 331.|18.}
\label{phi_\ind21_pic321b}
\end{center}
\end{figure}

\begin{minipage}[t]{0.5\textwidth}

\begin{enumerate}[start=332,leftmargin=15pt,label*={\arabic*}]

             \item [332.] If $(\vec{\w}\ |\vec{q}\ )\notin M_{\ind,2}^U$, then 
\\
\\

$$\cph\phi_{\ind|2,1}(\vec{\w}\ |\vec{q}\ ;z)=( \phi_{\ind|2}(\vec{\w}\ |\vec{q}\ ) ;z)\in M^U_{\ind+2,1}$$

\end{enumerate}

\end{minipage}
\vrule \hspace{0.5cm} 
\begin{minipage}[t]{0.5\textwidth}

\begin{enumerate}[start=1,leftmargin=5pt,label*={\arabic*}]

\item [19.] If $(\vec{\wii}\ |\xix) \in M^U_{\ind+2,1}$,
  $\xix\prec\wii_{\ind+2}$ and $\eta=0$, what implies
  $$ (\xix\prec\wii_{\theta+1}\text{ and }\xix\prec\wii_{\theta}),$$
  then we have
$$\cps\psi^2_{\ind|2,1}(\vec{\wii}\ ;\xix)=(\psi_{\ind|2}(\vec{\wii}\ );\xix).$$

\end{enumerate}

\end{minipage}
\bigskip

\begin{figure}[H]
\begin{center}
\begin{tikzpicture} [>=latex,every node/.style={minimum width=3em, node distance=4em},scale=1.5,transform shape]

\node [draw,circle, fill=blue!20] (b) {$\Scale[0.6]{\w_1\mid\wii_1}$};
\node [node distance=2.5em,right of=b] (d) {$...$}; 
\node [draw,circle, fill=orange!20, node distance=2.5em,right of=d] (f){$\Scale[0.6]{q_{01}\mid\wii_{\breve{\theta}}}$};
\node [node distance=2.5em,rotate=-54, above of=f, rotate=90] (mm1) {$\succ$}; 
\node [draw,circle, fill=orange!20, node distance=5em,right of=f] (fff) {$\Scale[0.6]{q_{01}\mid\wii_{\breve{\theta}+1}}$};
\node [node distance=2.5em,rotate=-56, above of=fff, rotate=90] (mm2) {$\succ$}; 
\node [node distance=2.5em,right of=fff] (g) {$\prec$}; 
\node [node distance=1em,right of=g] (h) {$...$}; 
\node [node distance=1em,right of=h] (i) {$\prec$}; 
\node [draw,circle, fill=blue!20, node distance=2.5em,right of=i] (ii) {$\Scale[0.6]{\w_{\ind-1}\mid\wii_{\ind+1}}$}; 
\node [node distance=2.5em,right of=ii] (iii) {$\prec$}; 
\node [draw,circle, fill=blue!20, node distance=2.5em,right of=iii] (j) {$\Scale[0.6]{\w_{\ind}\mid\wii_{\ind+2}}$};  
\node [node distance=2.5em, above of=j,rotate=-90] (k) {$\prec$}; 
\node [draw,circle, fill=orange!50, node distance=2.5em, above of=k] (l)  {$\Scale[0.6]{z\mid\xix}$}; 
\draw [<->] (f) -- (fff);

\end{tikzpicture}
\captionsetup{justification=centering,margin=2cm}
\caption{Illustration of 332.|19. (general case)}
\label{psi1_\ind21_pic6}
\end{center}
\end{figure}
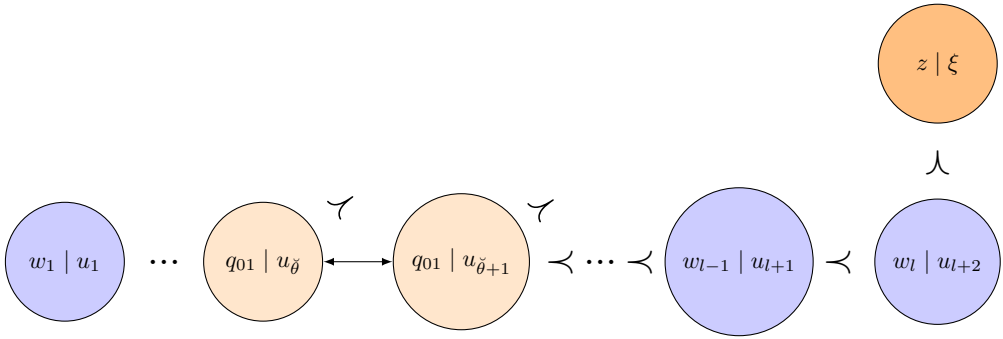

It is easy to see that by construction the left hand side fully describes the set $P^U_{\ind}\times M^U_{2,1}$.
Now, we check that the right hand side coincide with $M_{\ind+2,1}^U\sqcup M_{\ind+1,2}^U\sqcup M^U_{\ind,2,1}$:
\begin{itemize}
\item Let $(\vec{\wii}\ ; \xix)\in M^U_{\ind+2,1}$. Then the following cases from the right hand side clearly describe $M^U_{\ind,2,1}$:
\begin{enumerate}[leftmargin=40pt]
\item [2.] $\xix\succ\vec{\wii};$
\item [4.] $\xix\prec\wii_{\ind+2}$ and $\eta>0$;
\item [7.] $\xix\prec\wii_{\ind+2}$, $\eta=0$ and $(\xix\succ\wii_{\theta+1} \text{ and }  \xix\succ\wii_{\theta})$;
\item [19.] $\xix\prec\wii_{\ind+2}$, $\eta=0$ and $(\xix\prec\wii_{\theta+1} \text{ and }  \xix\prec\wii_{\theta})$.
\end{enumerate}

\bigskip

\item Let $(\vec{\w}\ ,\xix;\kap_0,\kap_1)\in M^U_{\ind+1,2}$. Then the following cases from the right hand side clearly describe $M^U_{\ind,2,1}$.:
\begin{enumerate}[leftmargin=40pt]

\item [3.] ($\xix\succ w_\ind,$ $\w_{\theta}\succ \kap_0,$ $\w_{\theta+1}\succ\kap_1,$ and $\xix\succ \kap_0$), or ($\xix\sim w_\ind,$ $\w_{\ind}\succ \kap_0,$ $\xix\succ\kap_1,$ and $\xix\succ \kap_0$);

\item [11.] $(\xix\sim w_\ind,$ $\w_{\ind}\succ \kap_0,$ $\xix\succ\kap_1,$ and $\xix\sim \kap_0);$
\bigskip
\item [9.] $q_0\succ\vec{w}$, $q_1\succ\xix$ and $q_0\sim\xix;$ 
\item [16.] $q_0\succ\vec{w}$, $q_1\succ\xix$ and $q_0\succ\xix.$ 
\end{enumerate}

\bigskip

\item Let $(\vec{\w}\ ;\kap_0,\kap_1; z)\in M^U_{\ind,2,1}$. Then the following cases from the right hand side clearly describe $M^U_{\ind,2,1}$.:
\begin{enumerate}[leftmargin=40pt]

\item [1.] $z\succ\vec{q}$ and $z\succ{w}$;
\item [5.] $z\succ\vec{q}$ and $z\prec{w_\ind}$;
\bigskip
\item [8.] $q_0\sim z$, $q_1\succ z$ and $z\succ\vec{w}$;
\item [12.] $q_0\sim z$, $q_1\succ z$ and $(w_\theta\succ q_0\text{ and }w_{\theta+1}\succ q_1)$;
\bigskip
\item [15.] $q_0\succ z$, $q_1\succ z$ and $z\prec w_\ind$;
\item [18.] $q_0\succ z$, $q_1\succ z$ and $z\prec w_\ind$;
\bigskip
\item [6.] First exceptional type element;
\item [13.] Second exceptional type element.
\end{enumerate}
\end{itemize}
Since every number from 1 to 19 was used exactly once, this completes the proof.

\end{proof}

\bigskip

\begin{thm}\label{2\ind1k}
For natural numbers $\ind$ and $k$, let $$M^U_{2^\ind,1^k}=\{(\vec{\xi},\vec{\varepsilon}\ )\in E_\ind^U\times E_{k+\ind}^U|\; \exists\{i_j\}_{j=1}^\ind,  \ 0<i_1<i_2<...<i_\ind<k+\ind, \text{ s.t. } \xi_j\sim\varepsilon_{i_j} \text{ for } 1\leq j\leq\ind  \},$$
Then
$$m_{2^\ind,1^k}^U=\sum_{(\vec{\xi},\vec{\varepsilon}\,)\in M^U_{2^\ind,1^k}}\xi_1\cdot...\cdot \xi_{\ind}\cdot \varepsilon_1\cdot...\cdot\varepsilon_{\ind+k}.$$ 

\end{thm}
\begin{remark}
According to Remark~\ref{c_m}, this implies $c_{2^k,1^{n-2k}}(U)\geq 0$.

\end{remark}
The proof is omitted and can be found in~\cite{Paunov16}.


\end{document}